\input ./style/arxiv-general.cfg
\documentclass[aos,MSNbibl,seceqn,dvips]{arximspdf}
\makeatletter
   \@ifpackageloaded{graphicx}{}{\usepackage{graphicx}}
\makeatother
\usepackage{dcolumn}

\doi{10.1214/15-AOS1313}
\volume{43}
\issue{3}
\pubyear{2015}
\firstpage{1351}
\lastpage{1390}
\docsubty{FLA}

\makeatletter
\newcommand{\Rem}{\mathrm{Rem}}
\newcommand{\Var}{\operatorname{Var}}

\newcommand{\lleft}{\left}
\newcommand{\rright}{\right}
\newcolumntype{d}[1]{D{.}{.}{#1}}
\newtheorem{theorem}{Theorem}[section]
\newtheorem{Lemma}[Theorem]{Lemma}
\newtheorem{Corollary}[Theorem]{Corollary}
\newtheorem{Proposition}[Theorem]{Proposition}
\newtheorem{Theorem}{Theorem}[section]
\newproclaim{Assumption}{Assumption}
\newtheorem{lemma}{Lemma}[section]
\newproclaim{Remark}{Remark}[section]
\newproclaim{Example}{Example}[section]
\makeatother
\begin{document}
\begin{frontmatter}

\title{Joint asymptotics for semi-nonparametric regression models with
partially linear structure}
\runtitle{Joint asymptotics for semi-nonparametric models}

\begin{aug}
\author[A]{\fnms{Guang}~\snm{Cheng}\thanksref{T1}\corref{}\ead[label=e1]{chengg@purdue.edu}}
\and
\author[A]{\fnms{Zuofeng}~\snm{Shang}\ead[label=e2]{shang9@purdue.edu}}
\runauthor{G. Cheng and Z. Shang}
\affiliation{Purdue University}
\address[A]{Department of Statistics\\
Purdue University\\
250 N. University Street\\
West Lafayette, Indiana 47907\\
USA\\
\printead{e1}\\
\phantom{E-mail:\ }\printead*{e2}}
\end{aug}
\thankstext{T1}{Supported by NSF CAREER Award DMS-11-51692,
DMS-14-18042 and Simons
Foundation Grant 305266.}

%
\received{\smonth{5} \syear{2014}}
%
\revised{\smonth{1} \syear{2015}}

%
\begin{abstract}
We consider a joint asymptotic framework for studying
semi-nonpara\-metric regression models where (finite-dimensional)
Euclidean parameters and (infinite-dimensional) functional parameters
are both of interest. The class of models in consideration share a
partially linear structure and are estimated in two general contexts:
(i) quasi-likelihood and (ii) true likelihood. We first show that the
Euclidean estimator and (pointwise) functional estimator, which are
re-scaled at different rates, jointly converge to a zero-mean Gaussian
vector. This weak convergence result reveals a surprising \textit{joint
asymptotics phenomenon}: these two estimators are asymptotically
independent. A~major goal of this paper is to gain first-hand insights
into the above phenomenon. Moreover, a likelihood ratio testing is
proposed for a set of joint local hypotheses, where a new version of
the Wilks phenomenon [\textit{Ann. Math. Stat.} \textbf{9} (1938)
60--62; \textit{Ann. Statist.} \textbf{1} (2001) 153--193] is
unveiled. A~novel
technical tool, called a \textit{joint Bahadur representation}, is
developed for studying these joint asymptotics results.
\end{abstract}

%
\end{frontmatter}

\section{Introduction}
In the literature, a statistical model is called \textit
{semi-nonpara\-metric} if it contains both finite-dimensional and
infinite-dimensional unknown parameters of interest (e.g., \cite
{HL07}). An example is semi-nonparametric copula model that can be
applied to address tail dependence among shocks to different financial
series and also to recover the shape of the ``news impact curve'' for
individual financial series. Another example is the semi-nonparametric
binary regression models proposed by Banerjee, Mukherjee and Mishra
\cite{BMM09} to define the conditional probability of attending primary
school in Indian villages through an appropriate link function
influenced by a set of covariates such as gender and household income.
As a first step in exploring the joint asymptotics results, we focus on
the semi-nonparametric regression models with a partial linear
structure in this paper.

The existing semiparametric literature is concerned with asymptotic
theories and inference procedures for the Euclidean parameter \textit
{only}. The functional parameter is profiled out as an
infinite-dimensional nuisance parameter; see \cite
{BKRW98,MV00,CH10,CK08,CK09,SS94}. In the special case where both
parameters are
estimable at the same root-n rate (e.g., \cite{KLF04,LPL08}), we can
combine them as an infinite-dimensional parameter and then apply the
functional Z-estimation theorem (e.g., Theorem~3.3.1 in \cite{VW96}),
to study its joint asymptotic distribution. However, it is more common
for the two parameters to be estimated at different parametric and
nonparametric rates. In general, their radically different parameter
dimensionality poses technical challenges for the construction of valid
procedures for joint inference. In this paper, we develop a new
technical tool, called a \textit{joint Bahadur representation} (JBR), for
studying the joint asymptotics results. As far as we are aware, our
joint asymptotic theories and inference procedures are new. The only
relevant reference of which we are aware is \cite{R08}, which focuses
on a fully parametric setting.

In this paper, we assume a partially linear structure for the
conditional mean of the response, and then estimate the model in two
general contexts: (i) quasi-likelihood and (ii) true likelihood. Within
this framework, we derive a joint limit distribution for the Euclidean
estimator and the (point-wise) functional estimator as a zero-mean
Gaussian vector after they are re-scaled properly. One surprising
result is that these two estimators are asymptotically independent.
This asymptotic independence will prove to be useful in making joint
inference. For example, it is now straightforward to construct the
joint confidence interval based on two marginal ones. Under similar
conditions, the marginal limit distribution for the Euclidean estimator
coincides with that derived in \cite{MVG97}. On the other hand, we
observe that the (pointwise) marginal asymptotic results for the
nonparametric component are generally different from those derived in a
purely nonparametric setup (without the Euclidean parameter) (i.e.,
\cite{SC12}), even though the Euclidean parameter is estimated at a
faster rate; see Remark~\ref{remeg}. This conclusion is a bit counterintuitive.

We next propose likelihood ratio testing for a variety of joint local
hypotheses such as $H_0\dvtx  \theta=\theta_0$ and $g(z_0)=w_0$ and $H_0\dvtx
x^T\theta+g(z_0)=\alpha$, where $\theta$ and $g$ denote the parametric
and nonparametric components, respectively. Conventional semiparametric
testing only focuses on the parametric components; see \cite
{MV00,CK09}. However, in practice, it is often of great interest to evaluate
the nonparametric components at the same time. For example, we may test
the joint effect of child gender $\theta$ and household income $g$ on
the probability of attending primary school in the Indian schooling
model; see \cite{BMM09}. In particular, we show that the null limit
distribution is a mixture of two \textit{independent} Chi-square
distributions that are contributed by the parametric and nonparametric
components, respectively. Note that this independence property is
implied by the joint asymptotics phenomenon, and is practically useful
in finding the critical value. In the parametric framework, Wilks
(1938) showed that the likelihood ratio test statistic (under $H_0\dvtx
\theta=\theta_0$) converges to a Chi-square distribution. Fan et al.
(2001) call the above result the \textit{Wilks phenomenon} due to the nice
property that the asymptotic null distribution is free of nuisance
parameters, and further generalize it to the nonparametric setting.
Therefore, we unveil a new version of Wilks phenomenon that adapts to
the semi-nonparametric context in this paper. As far as we are aware,
this joint testing result is new. The only relevant paper of which we
are aware is \cite{BMM09}, where the authors consider two separate null
hypotheses, that is, $H_{01}\dvtx  \theta=\theta_0$ and $H_{02}\dvtx  g(z_0)=w_0$,
under the monotonicity constraint of $g(\cdot)$.

The class of semi-nonparametric regression models considered in this
paper serves as a natural platform to deliver a new theoretical
insight: joint asymptotics phenomenon. We also note that our results
may be extended to the other models: (i) generalized additive partially
linear models, (ii) partial functional linear regression models \cite
{S09} and (iii) partially linear Cox proportional hazard models \cite
{H99} by either modifying the JBR or the criterion function; see
Section~\ref{secfut} for more elaborations. All the possible
extensions mentioned above require a smoothness assumption on the
nonparametric function. This assumption is crucially different from the
shape-constraint assumption, which in general leads to the
``nonstandard asymptotics'' problems (e.g., \cite{LN13,KP90,C09}).
Our framework cannot be easily adapted to handle these challenging
problems, which are usually analyzed by rather different technical tools.

The rest of this paper is organized as follows. Section~\ref{prel}
introduces the model assumptions and builds a theoretical foundation.
Sections~\ref{secthry} and~\ref{seclrt} formally discuss the joint
limit distribution and joint local hypothesis testing, respectively. In
Section~\ref{secexa}, we give three concrete examples with extensive
simulations to illustrate our theory. Section~\ref{secfut} discusses
some possible extensions. The proofs are postponed to the \hyperref[app]{Appendix} or
online supplementary document \cite{CS15}.

\section{Preliminaries}\label{prel}
This section introduces the model assumptions and establishes the
theoretical foundation of our results in two layers: (i) the partially
linear extension of reproducing kernel Hilbert space (RKHS) theory and
(ii) the joint Bahadur representation. Both technical results are of
independent interest.

\subsection{Notation and model assumptions}
Suppose that $T_i=(Y_i,X_i,Z_i)$, $i=1,\ldots,n$, are
i.i.d. copies of $T=(Y,X,Z)$, where
$Y\in\mathcal{Y}\subseteq\mathbb{R}$ is the response variable,
$U=(X,Z)\in\mathcal{U}\equiv\mathbb{I}^p\times\mathbb{I}$ is the
covariate variable, and $\mathbb{I}=[0,1]$.
Throughout the paper we assume
that the density of $Z$, denoted by $\pi(z)$,
has positive lower bound and finite upper bound for $z\in[0,1]$.
Consider a general class of semi-nonparametric
regression models with the following partially linear structure:
\begin{equation}
\label{partlinearmodel} \mu_0(U)\equiv E(Y|U)=F \bigl(X^T
\theta_0+g_0(Z) \bigr),
\end{equation}
where $F(\cdot)$ is some known link function and $g_0(\cdot)$ is some
unknown smooth function. This primary assumption covers two classes of
statistical models. The first class is called \textit{generalized
partially linear models} \cite{BHZ06}; here the data are modeled by
$y|u\sim p(y;\mu_0(u))$ for a conditional distribution $p$. Instead of
assuming the underlying distribution, the second class specifies only
the relationship between the conditional mean and the conditional
variance: $\operatorname{Var}(Y|U)=\mathcal{V}(\mu_0(U))$ for some known
positive-valued function $\mathcal{V}$. The nonparametric estimation of
$g$ in the second situation uses the quasi-likelihood $Q(y; \mu)\equiv
\int_{y}^\mu(y-s)/\mathcal{V}(s)\,ds$ with $\mu=F(x^T\theta+g(z))$
\cite{W74}. Despite the distinct modeling principles, these two classes have
a large overlap under many common combinations of $(F, \mathcal V)$, as
summarized in Table~2.1 of \cite{MN89}. From now on, we work with a
general criterion function $\ell(y;a)\dvtx \mathcal Y\times\mathbb
R\mapsto
\mathbb R$, which can represent either $\log p(y; F(a))$ or $Q(y;F(a))$.

Let the full parameter space for $f\equiv(\theta,g)$ be $\mathcal
{H}\equiv\mathbb{R}^p\times H^m(\mathbb{I})$, where $H^m(\mathbb
{I})$ is an
$m$th order Sobolev space defined as
\begin{eqnarray*}
H^m(\mathbb{I})&\equiv& \bigl\{g \dvtx\mathbb{I} \mapsto \mathbb{R}|
g^{(j)} \mbox{ is absolutely continuous}
\\
&&{}\hspace*{5pt}\mbox{for $j=0,1,\ldots ,m-1$ and }g^{(m)} \in L_{2}(\mathbb{I}) \bigr\}.
\end{eqnarray*}
With some abuse of notation, $\mathcal H$
may also refer to $\mathbb{R}^p\times H^m_0(\mathbb{I})$, where
$H_0^m(\mathbb{I})$ is a homogeneous subspace of $H^m(\mathbb{I})$. The
space $H^m_0(\mathbb{I})$ is also known as the class of periodic
functions such that a function $g\in H_0^m(\mathbb{I})$ has additional
restrictions $g^{(j)}(0)=g^{(j)}(1)$ for $j=0,1,\ldots,m-1$.
Throughout this paper we assume $m> 1/2$ to be known. Consider the
penalized semi-nonparametric estimator
\begin{eqnarray}
 (\widehat{\theta}_{n,\lambda},\widehat{g}_{n,\lambda})&=&
\mathop{\arg\max}_{(\theta,g)\in\mathcal{H}}\ell_{n,\lambda}(f)
\nonumber
\\[-8pt]
\label{estf}
\\[-8pt]
\nonumber
&=&\mathop{\arg\max}_{(\theta,g)\in\mathcal{H}} \Biggl\{ \frac{1}{n} \sum
_{i=1}^n \ell \bigl(Y_i;
X_i^T\theta+g(Z_i) \bigr)-(\lambda/2)
J(g,g) \Biggr\},
\end{eqnarray}
where $J(g,\tilde{g})=\int_\mathbb{I} g^{(m)}(z)\tilde
g^{(m)}(z)\, dz$ and $\lambda\rightarrow0$ as $n\rightarrow\infty$.
Here, we use $\lambda/2$ (rather than $\lambda$) to simplify future
expressions. Write $\widehat{f}_{n,\lambda}=(\widehat{\theta
}_{n,\lambda
},\widehat{g}_{n,\lambda})$. The existence of $\widehat g_{n,\lambda}$
is guaranteed by Theorem~2.9 of \cite{G02} when the null space
$\mathcal{N}_m\equiv\{g\in H^m(\mathbb{I})\dvtx  J(g,g)=0\}$ is finite-dimensional and $\ell(y;a)$ is concave and continuous w.r.t. $a$.

We next assume some basic model conditions. For simplicity, throughout
the paper we do not distinguish $f=(\theta,g)\in\mathcal{H}$ from its
associated function $f\in\mathcal{F}\equiv\{f(x,z)=x^T\theta+g(z)\dvtx
(\theta, g)\in\mathcal{H}, (x,z)\in\mathcal{U}\}$. Let $\mathcal{I}_0$ be
the range for the true function $f_0(x,z)\in\mathcal F$, that is, a compact
interval. Denote the first-, second- and third-order derivatives of
$\ell(y;a)$ (w.r.t. $a$) by $\dot\ell_a$, $\ddot{\ell}_a$ and
$\ell'''_a$.

\renewcommand{\theAssumption}{A\arabic{Assumption}}
\begin{Assumption}\label{A.1}
\textup{(a)} $\ell(y;a)$ is three times continuously
differentiable and concave w.r.t. $a$. There exists a bounded open
interval $\mathcal{I}\supset\mathcal{I}_0$ and positive
constants $C_0$ and $C_1$ s.t.
\begin{equation}
\label{A1aeq1} E \Bigl\{\exp \Bigl(\sup_{a\in\mathcal{I}}\bigl|\ddot{\ell
}_a(Y;a) \bigr|/C_0 \Bigr)\big|U \Bigr\}\le C_1
\qquad \mbox{a.s.}
\end{equation}
and
\begin{equation}
\label{A1aeq2} E \Bigl\{\exp \Bigl(\sup_{a\in\mathcal{I}} \bigl|
\ell'''_a(Y;a)\bigr|/C_0
\Bigr) \big| U \Bigr\}\le C_1\qquad \mbox{a.s.}
\end{equation}

(b) There exists a positive constant $C_2$
s.t. $C_2^{-1}\le I(U)\equiv
-E(\ddot{\ell}_a(Y;X^T\theta_0+g_0(Z))|U)\le C_2$, a.s.

(c) $\epsilon\equiv
\dot{\ell}_a(Y;X^T\theta_0+g_0(Z))$ satisfies
$E(\epsilon|U)=0$, $E(\epsilon^2|U)= I(U)$, a.s., and $E\{\epsilon
^4\}
<\infty$.
\end{Assumption}

A detailed discussion of the above model assumptions can be
found in \cite{SC12}.
In particular, Assumption~\ref{A.1}(a) is typically
used in semiparametric quasi-likelihood models; see \cite{MVG97}. Three
concrete examples showing the
validity of Assumption~\ref{A.1} are presented in Section~\ref{secexa}.

Hereinafter, if for positive
sequences $a_\mu$ and $b_\mu$ we have that $a_\mu/b_\mu$ tends to a
strictly positive constant, we write $a_\mu\asymp b_\mu$. If that
constant is one,
we write $a_\mu\sim b_\mu$. Let $\sum_\nu$ denote the sum
over $\nu\in\mathbb{N}=\{0,1,2,\ldots\}$ for convenience. Let the
sup-norm of
$g\in H^m(\mathbb{I})$ be $\|g\|_{\sup}=\sup_{z\in\mathbb{I}}|g(z)|$.
Let $\lambda^\ast$ be the optimal smoothing parameter;
$\lambda^\ast\asymp n^{-2m/(2m+1)}$. For simplicity, we write
$\lambda
^{1/(2m)}$ as $h$, and thus $h^\ast\asymp n^{-1/(2m+1)}$.

\subsection{A partially linear extension of RKHS theory}\label{subsection2.3}
In this section, we adapt the nonparametric RKHS framework to our
semi-nonparametric setup.

We define the inner product for $\mathcal H$ to be, for any $(\theta
,g),(\tilde{\theta},\tilde{g})\in\mathcal{H}$,
\begin{equation}
\label{innprodH}\quad \bigl\langle(\theta,g),(\tilde{\theta},\tilde{g}) \bigr\rangle
=E_U \bigl\{I(U) \bigl(X^T\theta+g(Z) \bigr)
\bigl(X^T\tilde{\theta}+\tilde{g}(Z) \bigr) \bigr\} +\lambda J(g,
\tilde{g}),
\end{equation}
and we define the norm to be $\|(\theta,g)\|^2=\langle(\theta
,g),(\theta,g)\rangle$.
The validity of such a norm is demonstrated in Section~S.1 of the
supplement document
\cite{CS15} (under Assumption~\ref{A.2} introduced later).
Under this norm, we will construct two linear operators, $R_u\in
\mathcal{H}$, for any $u\in\mathcal{U}$, and $P_\lambda\dvtx
\mathcal{H}\mapsto
\mathcal{H}$ satisfying
\begin{equation}
\label{RPKH} \langle R_u, f\rangle=x^T\theta+g(z)
\qquad \mbox{for any $u\in\mathcal U$ and $f\in\mathcal{H}$}
\end{equation}
and
\begin{equation}
\label{Penprop} \langle P_\lambda f,\tilde{f}\rangle=\lambda J(g,
\tilde{g})\qquad \mbox{for any $f=(\theta,g), \tilde {f}=(\tilde{\theta},
\tilde{g})\in\mathcal{H}$}.
\end{equation}
As will be seen, $R_u$ and $P_\lambda$ are two major building blocks of
this enlarged RKHS framework. In particular, Propositions~\ref{proru}
and~\ref{Penopr} show that these two operators are actually built upon
their nonparametric counterparts $K_z$ and $W_\lambda$ defined below.

Let $K(z_1,z_2)$ be a (symmetric) reproducing kernel of $H^m(\mathbb{I})$
endowed with the inner product $\langle g,\tilde{g}\rangle
_1=E_Z\{B(Z)g(Z)\tilde{g}(Z)\}+\lambda
J(g,\tilde{g})$ and norm $\|g\|_1^2=\langle g, g\rangle_1$, where
$B(Z)=E\{I(U)|Z\}$. Hence, $K_z(\cdot)\equiv K(z,\cdot)$ satisfies
$\langle K_z, g\rangle_1=g(z)$. We next specify a positive definite
self-adjoint operator $W_\lambda\dvtx  H^m(\mathbb{I})\mapsto H^m(\mathbb
{I})$ satisfying $\langle W_\lambda g,\tilde{g}\rangle_1=\lambda
J(g,\tilde{g})$ for any $g,\tilde{g}\in H^m(\mathbb{I})$.
The existence of such $W_\lambda$ is proved in Section~S.2 of the
supplement document \cite{CS15}.
Write
$V(g,\tilde{g})=E_Z\{B(Z)g(Z)\tilde{g}(Z)\}$. Hence, $\langle
g,\tilde{g}\rangle_1=V(g,\tilde{g})+\langle
W_\lambda g,\tilde{g}\rangle_1$, which implies
\begin{equation}
\label{sec2eq2} V(g,\tilde{g})= \bigl\langle (id-W_\lambda)g, \tilde{g}
\bigr\rangle_1,
\end{equation}
where $id$ denotes the identity operator. We next assume that there
exists a sequence of basis
functions in the space $H^m(\mathbb{I})$ that simultaneously
diagonalizes the bilinear forms $V$ and $J$. Such an eigensystem
assumption is
typical in the smoothing spline literature; see \cite{G02}.

\begin{Assumption} \label{A.3}
There exists a sequence of real-valued functions $h_\nu\in
H^m(\mathbb{I})$, $\nu\in\mathbb{N}$ satisfying
$\sup_{\nu\in\mathbb{N}}\|h_\nu\|_{\mathrm{sup}}<\infty$ and a
nondecreasing real sequence $\gamma_\nu\asymp\nu^{2m}$
such that $V(h_\mu,h_\nu)=\delta_{\mu\nu}$ and $J(h_\mu,h_\nu
)=\gamma
_\mu\delta_{\mu\nu}$ for any $\mu,\nu\in\mathbb{N}$, where
$\delta_{\mu\nu}$
is the Kronecker's delta. Furthermore, any
$g\in H^m(\mathbb{I})$ admits the Fourier expansion
$g=\sum_{\nu} V(g,h_\nu)h_\nu$ under the $\|\cdot\|_1$-norm.
\end{Assumption}

Under Assumption~\ref{A.3} and by $B(Z)=E\{I(U)|Z\}$,
it can be seen that $E\{I(U)h_\nu(Z)h_\mu(Z)\}=V(h_\nu,h_\mu
)=\delta
_{\nu\mu}$.
Then we can easily derive explicit expressions for $\|g\|_1$,
$W_\lambda h_\nu(\cdot)$ and $K_z(\cdot)$ in terms of the $h_\nu$ and
$\gamma_\nu$ as follows:
\begin{eqnarray}
\hspace*{6pt}\quad\|g\|_1^2 &=& \sum_{\nu}
\bigl|V(g,h_\nu)\bigr|^2(1+\lambda \gamma_\nu),\qquad
W_\lambda h_\nu(\cdot)= \frac{\lambda\gamma_\nu}{1+\lambda\gamma_\nu}h_\nu(
\cdot )\quad\mbox{and}
\nonumber
\\[-8pt]
\label{PropRPK}
\\[-8pt]
\nonumber
K_z(\cdot) &=& \sum_\nu
\frac{h_\nu(z)}{1+\lambda\gamma_\nu} h_\nu(\cdot).
\end{eqnarray}
Using similar arguments to those in Proposition~2.2 of \cite{SC12}, we
know that Assumption~\ref{A.3} holds when Assumption~\ref{A.1} is
satisfied and the $h_\nu$s are chosen as the (normalized) solutions of
the following ODE problem:
\begin{eqnarray}
(-1)^m h_\nu^{(2m)}(\cdot) &=&
\gamma_\nu B(\cdot)\pi(\cdot) h_\nu(\cdot),
\nonumber
\\[-8pt]
\label{eigenproblem}
\\[-8pt]
\nonumber
h_\nu^{(j)}(0) &=& h_\nu^{(j)}(1)=0,
\qquad j=m,m+1,\ldots,2m-1.
\end{eqnarray}
For example, the $h_\nu$s are constructed as an explicit trigonometric
basis in case (I) of Example~\ref{egpss}.
As will be seen later, by employing the above ordinary differential
equation (ODE) approach, we will reduce the challenging
infinite-dimensional inference problems to simple exercises on finding
the underlying eigensystem. We remark that proving the existence of the
above eigensystem is nontrivial and relies substantially on the ODE
techniques developed in \cite{Birk1908,Stone1926}.

We next state a regularity Assumption \ref{A.2} guaranteeing that $R_u$
and $P_\lambda$ are both well defined. Define $A_0(Z)=E\{I(U)X|Z\}$ and
$G(Z)=A_0(Z)/B(Z)$. Note that $G=(G_1,\ldots,G_p)^T$ is a
$p$-dimensional vector-valued function, for example, $G(Z)=E(X|Z)$ in
the $L_2$ regression.

\begin{Assumption}\label{A.2}
$G_1,\ldots,G_p\in L_{2}(P_Z)$, that is,
$G_k$ has a finite second moment, and the $p\times p$ matrix $\Omega
\equiv
E\{I(U)(X-G(Z))(X-G(Z))^T\}$ is positive definite.
\end{Assumption}

Under the assumption that $G_k\in L_{2}(P_Z)$, the
linear functional $\mathcal{A}_k$ defined by
$\mathcal{A}_kg=V(G_k,g)$ is bounded (or equivalently, continuous) for
any $g\in
H^m(\mathbb{I})$ because of the following inequality:
$|\mathcal{A}_kg|\leq V^{1/2}(G_k, G_k)V^{1/2}(g,g)\leq V^{1/2}(G_k,
G_k)\|g\|_1<\infty$. Thus, by Riesz's representation
theorem, there exists an $A_k\in H^m(\mathbb{I})$ such that
$\mathcal{A}_kg=\langle A_k,g\rangle_1$ for any $g\in
H^m(\mathbb{I})$. Thus if we write $A=(A_1,\ldots,A_p)^T$, then
\begin{equation}
\label{inter1} V(G,g)=\langle A,g\rangle_1.
\end{equation}
We also note that $A=(id-W_\lambda)G$ when $G_1,\ldots,G_p$
$\in H^m(\mathbb{I})$ based on (\ref{sec2eq2}). Taking $g=K_z$ in
(\ref{inter1}) and applying (\ref{PropRPK}), we find that
\begin{equation}
\label{wlamA}\quad\hspace*{20pt} A(z)=\sum_{\nu}\frac{V(G,h_\nu)}{1+\lambda\gamma_\nu}h_\nu(z)
\quad\mbox{and}\quad (W_\lambda A) (z)=\sum_{\nu}
\frac{V(G,h_\nu)\lambda\gamma
_\nu
}{(1+\lambda\gamma_\nu)^2}h_\nu(z).\hspace*{-4pt}
\end{equation}

Now,\vspace*{1.5pt} we are ready to construct $R_u$ and $P_\lambda$ in
Propositions~\ref{proru} and~\ref{Penopr}, respectively. Define
$\Sigma_\lambda=E_Z\{B(Z)G(Z)(G(Z)-A(Z))^T\}$ as a $p\times p$ matrix.

\begin{Proposition}\label{proru}
$R_u$ defined in (\ref{RPKH}) can be expressed as $R_u\dvtx  u\mapsto
(H_u,T_u)\in\mathcal H$, where
\begin{eqnarray}
H_u &=& (\Omega+\Sigma_\lambda)^{-1}
\bigl(x-A(z) \bigr) \quad\mbox{and}
\nonumber
\\[-8pt]
\label{RTu}
\\[-8pt]
\nonumber
 T_u &=& K_z-A^T(
\Omega+\Sigma_\lambda)^{-1} \bigl(x-A(z) \bigr).
\end{eqnarray}
\end{Proposition}

\begin{Proposition}\label{Penopr}
$P_\lambda$ defined in (\ref{Penprop}) can be expressed as
$P_\lambda\dvtx
(\theta, g)\mapsto(H_g^\ast, T_g^\ast)\in\mathcal H$, where
\[
\cases{ H_g^*=-(\Omega+\Sigma_\lambda)^{-1}E
\bigl\{ B(Z)G(Z) (W_\lambda g) (Z) \bigr\},\vspace*{3pt}
\cr
T_g^*=E \bigl\{B(Z)G(Z)^T(W_\lambda g) (Z) \bigr
\}(\Omega+\Sigma_\lambda )^{-1}A+W_\lambda g.}
\]
\end{Proposition}

Note that $\lim_{\lambda\rightarrow0}\Sigma_\lambda=0$
according to (\ref{lemma3ii}) in the \hyperref[app]{Appendix}. Therefore, $(\Omega
+\Sigma_\lambda)^{-1}$ above is well defined under Assumption~\ref
{A.2}. In addition, we note that $P_\lambda$ is self-adjoint and
bounded because of the following inequality:
\begin{eqnarray}
\label{boundPlambda} \|P_\lambda f\|&=&\sup_{\|\tilde{f}\|= 1} \bigl|\bigl\langle
P_\lambda f,\tilde{f}\bigr\rangle \bigr|
\nonumber
\\[-8pt]
\\[-8pt]
\nonumber
&=&\sup_{\|\tilde{f}\|= 1} \bigl|\lambda J(g,\tilde{g}) \bigr|\le\sqrt{
\lambda J(g,g)} \sup_{\|\tilde{f}\|=1} \sqrt{\lambda J(\tilde{g},
\tilde{g})}\le\|f\|.
\end{eqnarray}

Finally, we derive the Fr\'{e}chet derivatives of $\ell_{n,\lambda}(f)$
defined in (\ref{estf}). Let $\Delta f, \Delta f_j\in\mathcal H$
for $j=1,2,3$. The Fr\'{e}chet derivative of $\ell_{n,\lambda}(f)$ is
given by
\begin{eqnarray*}
D\ell_{n,\lambda}(f)\Delta f&=&\frac{1}{n}\sum
_{i=1}^{n}\dot{\ell}_a
\bigl(Y_i; X_i^T\theta+g(Z_i)
\bigr) \langle R_{U_i},\Delta f\rangle-\langle P_\lambda f,\Delta
f\rangle
\\
&\equiv& \bigl\langle S_n(f), \Delta f \bigr\rangle-\langle
P_\lambda f,\Delta f\rangle\equiv \bigl\langle S_{n,\lambda}(f),\Delta
f \bigr\rangle.
\end{eqnarray*}
Note that $S_{n,\lambda}(\widehat{f}_{n,\lambda})=0$. In particular,
$S_{n,\lambda}(f_0)$ is of interest, and it can be expressed as
\begin{equation}
\label{snlam} S_{n,\lambda}(f_0)=\frac{1}{n} \sum
_{i=1}^n\epsilon_i R_{U_i}-P_\lambda
f_0.
\end{equation}
The Frech\'{e}t derivatives of $S_{n,\lambda}$ and  $\mathit{DS}_{n,\lambda}$,
denoted $\mathit{DS}_{n,\lambda}(f)\Delta f_1\Delta f_2$ and $D^2S_{n,\lambda}(f)
\Delta f_1\Delta f_2 \Delta f_3$, can be explicitly calculated as
$(1/n)\sum_{i=1}^n
\ddot{\ell}_a(Y_i;\break X_i^T\theta+g(Z_i))\langle R_{U_i}, \Delta
f_1\rangle\langle R_{U_i},\Delta f_2\rangle-\langle P_\lambda\Delta
f_1,\Delta f_2\rangle$ and $(1/n)\sum_{i=1}^n
\ell'''_a(Y_i;\break X_i^T\theta+   g(Z_i))\langle R_{U_i}, \Delta
f_1\rangle\langle R_{U_i},\Delta f_2\rangle\langle R_{U_i},\Delta
f_3\rangle$, respectively. Define
$S(f)=E\{S_n(f)\}$, $S_\lambda(f)=S(f)-P_\lambda f$ and $\mathit{DS}_\lambda
(f)=\mathit{DS}(f)-P_\lambda$,\vspace*{1pt}
where\break $\mathit{DS}(f)\Delta f_1\Delta f_2=E\{\ddot{\ell}_a(Y;
X^T\theta+g(Z))\langle R_{U}, \Delta f_1\rangle\langle
R_{U},\Delta f_2\rangle\}$. Since\break $\langle \mathit{DS}_\lambda(f_0) f,
\tilde{f}\rangle=-\langle f,\tilde{f}\rangle$ for any
$f,\tilde{f}\in\mathcal{H}$, we have the following result:

\begin{Proposition}\label{basicprop}
$\mathit{DS}_\lambda(f_0)=-id$, where $\mathit{id}$ is the identity
operator on $\mathcal{H}$.
\end{Proposition}

\subsection{Joint Bahadur representation}\label{sectech}

This section presents the second layer of our theoretical foundation,
the \textit{joint Bahadur representation} (\textit{JBR}). The JBR is developed
based on empirical processes theory and will prove to be a powerful
tool in the study of joint asymptotics.

We start with a useful lemma stating the relationship between $\|f\|$
and $\|f\|_{\mathrm{sup}}$, where the former $f=(\theta, g)$ and the latter
$f=x^T\theta+g(z)$.

\begin{Lemma}\label{lemma0}
There exists a constant $c_m>0$ such that $\|R_u\|\leq c_mh^{-1/2}$ and
$|f(u)|\le
c_m h^{-1/2}\|f\|$ for any $u\in\mathcal{U}$ and
$(\theta,g)\in\mathcal{H}$. In particular, $c_m$ does not
depend on the choice of $u$ and $(\theta,g)$. Hence
$\|f\|_{\mathrm{sup}}\leq c_m h^{-1/2}\|f\|$.
\end{Lemma}

An additional\vspace*{1pt} convergence-rate condition is needed to obtain JBR.
Assumption~\ref{A.5} implies that $\widehat{f}_{n,\lambda}$ achieves the
optimal rate of convergence, that is, $O_P(n^{-m/(2m+1)})$, when
$\lambda=\lambda^*$.

\begin{Assumption}\label{A.5}
$
\|\widehat{f}_{n,\lambda}-f_0\|=O_{P}
((nh)^{-1/2}+h^m).
$
\end{Assumption}

\noindent Interestingly, we show below that the above rate condition
(Assumption~\ref{A.5}) is valid for a broad range of $h$ once
Assumptions \ref{A.1}--\ref{A.2} hold (by employing the contraction
mapping idea).

\begin{Proposition}\label{ratesconvergence}
Suppose Assumptions \textup{\ref{A.1}}--\textup{\ref{A.2}} are satisfied. Furthermore, as
$n\rightarrow\infty$,
$h=o(1)$ and $n^{-1/2}h^{-2}(\log{n})(\log\log{n})^{1/2}=o(1)$.
Then $\|\widehat{f}_{n,\lambda}-f_0\|=O_P((nh)^{-1/2}+h^m)$.
\end{Proposition}

We remark that the optimal rate for the smoothing parameter,
that is, $h^\ast\asymp n^{-1/(2m+1)}$, satisfies the rate conditions
for $h$ specified in Proposition~\ref{ratesconvergence} when $m>3/2$.

The following \textit{joint Bahadur representation} can be viewed as a
nontrivial extension of the Bahadur representation \cite{B66} for
parametric models by adding a functional component.

\begin{Theorem}[(Joint Bahadur representation)]\label{mainthm1}
Suppose that Assumptions \ref{A.1} through \ref{A.5} hold, $h=o(1)$ and
$nh^2\rightarrow\infty$. Recall that $S_{n,\lambda}(f_0)$ is defined
in (\ref{snlam}). Then we have
\begin{equation}
\label{bahareps}
\bigl\|\widehat{f}_{n,\lambda}-f_0-S_{n,\lambda}(f_0)
\bigr\|=O_{P}(a_n\log{n}),
\end{equation}
where $a_n=n^{-1/2}((nh)^{-1/2}+h^m)
h^{-(6m-1)/(4m)}(\log\log{n})^{1/2}+\break C_\ell h^{-1/2}
((nh)^{-1}{+} h^{2m})/\log{n}$ and
$C_\ell=\sup_{u\in\mathcal{U}}E\{\sup_{a\in\mathcal{I}}|\ell
'''_a(Y;a)| |U=u\}$.
\end{Theorem}

The proof of Theorem~\ref{mainthm1} relies heavily on modern
empirical process theory, and in particular a \textit{concentration
inequality} given in the supplementary material~\cite{CS15}.

%

\section{Joint limit distribution}\label{secthry}

As far as we are aware, the current semiparametric literature on the
smoothing spline models mostly focus on the asymptotic normality of the
parametric parts, and derive only rates of convergence (in estimation)
for functional parts; see \cite{W90,G02,W11}. In this section, we
demonstrate the joint asymptotic normality of both parametric and
functional parts.

We start from a preliminary result that for any $z_0\in\mathbb{I}$,
$(\sqrt{n}(\widehat{\theta}_{n,\lambda}-\theta_0^\ast),\break \sqrt
{nh}(\widehat{g}_{n,\lambda}-g_0^\ast)
(z_0))$ weakly converges to a zero-mean Gaussian vector. Unfortunately,
the center $(\theta_0^*,g_0^*)\equiv f_0-P_{\lambda}f_0$ is biased and
the asymptotic variance is not diagonal; see Theorem~\ref{mainthm2} in
the \hyperref[app]{Appendix} for more technical details. Under a regularity condition
on the least favorable direction \cite{K08}, that is, (\ref{smoothG}),
we can remove the estimation bias for $\theta$; see Lemma~\ref{lemma4}
in the \hyperref[app]{Appendix}. In this case, the parametric estimate $\widehat\theta
_{n,\lambda}$ is semiparametric efficient when $Y$ belongs to an
exponential family; see Remark~\ref{remcom1}. However, what is more
surprising is that \textit{$\widehat\theta_{n,\lambda}$ and $\widehat
g_{n,\lambda}(z_0)$ become asymptotically independent} after the bias
removal procedure. We call this discovery the \textit{joint asymptotics
phenomenon}. This leads to the first main result of this paper, given
in Theorem~\ref{mainthm3} below.

\begin{Theorem}[(Joint limit distribution)]\label{mainthm3}
Let Assumptions \ref{A.1} through \ref{A.5} be satisfied.
Suppose there exists $b\in(1/(2m),1]$ such that $G_k$ satisfies
\begin{equation}
\label{smoothG} \sum_{\nu}\bigl|V(G_k,h_\nu)\bigr|^2
\gamma_\nu^b<\infty\qquad \mbox{for any $k=1,\ldots,p$}.
\end{equation}
Furthermore,\vspace*{1pt} as $n\rightarrow\infty$, $h=o(1)$,
$nh^2\rightarrow\infty$, $a_n\log{n}=o(n^{-1/2}h^{1/2})$
[with $a_n$ defined as in (\ref{bahareps})],
$h V(K_{z_0},K_{z_0})\rightarrow\sigma^2_{z_0}$ and
$n^{1/2}h^{m(1+b)}=o(1)$.
Then we have, for any $z_0\in\mathbb I$,
\begin{equation}
\label{importantlimit}
\pmatrix{\sqrt{n}(\widehat{
\theta}_{n,\lambda}-\theta_0)
\vspace*{2pt}\cr
\sqrt{nh} \bigl\{\widehat{g}_{n,\lambda}(z_0)-g_0(z_0)+(W_\lambda
g_0) (z_0) \bigr\} }
\stackrel{d} {\longrightarrow} N (0,\Psi ),
\end{equation}
where
\begin{equation}
\label{Psieqmainthm3}
\Psi=\pmatrix{
\Omega^{-1} & 0
\vspace*{2pt}\cr
0& \sigma_{z_0}^2 }.
\end{equation}
Furthermore, if
\begin{equation}\label{weta}
\lim_{n\rightarrow\infty}(nh)^{1/2}(W_\lambda
g_0) (z_0)=-b_{z_0},
\end{equation}
then we have
\begin{equation}
\label{animportantlimit}
\pmatrix{
\sqrt{n}(\widehat{
\theta}_{n,\lambda}-\theta_0)
\vspace*{2pt}\cr
\sqrt{nh} \bigl(\widehat{g}_{n,\lambda}(z_0)-g_0(z_0)
\bigr) }
  \stackrel{d} {\longrightarrow} N
\biggl(\pmatrix{0 \vspace*{0pt}\cr b_{z_0}},\Psi \biggr).
\end{equation}
\end{Theorem}

We remark that Theorem~\ref{mainthm3} holds under the
optimal smoothing parameter $h^\ast=n^{-1/(2m+1)}$. It follows from
(\ref{PropRPK}) and (\ref{wlamA}) that
\begin{equation}\label{sigma}
\sigma_{z_0}^2=\lim_{h\rightarrow0}\sum
_{\nu}\frac
{h|h_\nu(z_0)|^2}{(1+\lambda\gamma_\nu)^2}.
\end{equation}

It is worth pointing out that we obtain the above results without
strengthening the regularity conditions used in the semiparametric
literature, for example, those in \cite{MVG97}. We next discuss the key
condition (\ref{smoothG}).
When $b=0$, condition (\ref{smoothG}) reduces to Assumption~\ref{A.2}
that $G_k\in L_{2}(P_Z)$. However, we require
$1/(2m)<b\le1$ such that the Fourier coefficients $V(G_k,h_\nu)$ in
(\ref{smoothG}) converge to zero at a faster rate than $\nu^{-mb}$
because $\gamma_\nu\asymp\nu^{2m}$; see Assumption~\ref{A.3}. It is
well known that a faster decaying rate of
the Fourier coefficients $V(G_k,h_\nu)$ implies a more smooth $G_k$;
see \cite{CO90}, page 1681.
Therefore, condition (\ref{smoothG}) requires more
smoothness of $G_k$. In fact, it follows from \cite{CO90} that (\ref
{smoothG}) is
equivalent to $G_k\in H^{mb}(\mathbb{I})$ with $1/(2m)<b\le1$. Hence,
the condition $G_k\in H^{m}(\mathbb{I})$
assumed in the classical semiparametric work by Mammen and van de Geer
\cite{MVG97} may actually be weakened.

We next discuss three important consequences of Theorem~\ref
{mainthm3}. First, the asymptotic independence between $\widehat\theta
_{n,\lambda}$ and $\widehat g_{n,\lambda}(z_0)$ greatly facilitates the
construction of the joint CI for $(\theta_0, g_0(z_0))$ by directly
building on the marginal CIs. Second, based on Theorem~\ref{mainthm3}
and the Delta method, we can easily establish the prediction interval
for a new response $Y_{\mathrm{new}}$ given future data $u_0=(x_0,z_0)$ and the
CI for some real-valued smooth function of $(\theta_0, g_0(z_0))$; see
Section~\ref{secexa}.
Finally, the nonparametric estimation bias $b_{z_0}$ can be further
removed under an additional assumption; see Corollary~\ref{l2regressionCLT}.

In Remarks~\ref{remcom1} and~\ref{remcom2} below, we compare the
marginal limit distributions implied by Theorem~\ref{mainthm3} with
those derived in the semiparametric \cite{MVG97} and nonparametric
\cite
{SC12} literature.

\begin{Remark}\label{remcom1}
Our parametric limit distribution is $\sqrt{n}(\widehat\theta
_{n,\lambda
}-\theta_0)\stackrel{d}{\longrightarrow}N(0,
\Omega^{-1})$, where $\Omega=E\{I(U)(X-G(Z))(X-G(Z))^T\}$. We find that
it is exactly the same as that obtained in \cite{MVG97}; see
Section~S.14 of supplementary document \cite{CS15}. Mammen and van de
Geer \cite{MVG97} further showed that the parametric estimate is
semiparametric efficient when $Y$ belongs to an exponential family; see
their Remark~4.1. For example, in the partially linear models under
Gaussian errors, $\Omega$~reduces to the semiparametric efficiency
bound $E(X-E(X|Z))^{\otimes2}$; see \cite{K08}. Note the profile
approach in \cite{MVG97} treats $g$ as a nuisance parameter, and thus
it cannot be adapted to obtain our joint limiting distribution.
\end{Remark}

\begin{Remark}\label{remcom2}
Our\vspace*{1pt} (pointwise) nonparametric limit distribution, that is, $\sqrt
{nh}(\widehat{g}_{n,\lambda}(z_0)-g_0(z_0))\stackrel
{d}{\longrightarrow
}N(b_{z_0}, \sigma^2_{z_0})$, is in general different from that
obtained in the nonparametric smoothing spline setup (without $\theta$)
in terms of different values of $b_{z_0}$ and $\sigma_{z_0}^2$; see
\cite{SC12}. This is mainly due to the eigensystem difference in the
two setups; see Remark~\ref{remeg} for more illustrations. An
exception is the $L_2$ regression in which the two eigensystems
coincide. Our general finding gives a counter-example to the common
intuition in the literature that the nonparametric limit distribution
is not affected by the involvement of a parametric component that is
estimated at a faster convergence rate.
\end{Remark}

To further illustrate Theorem~\ref{mainthm3}, we consider the partial
smoothing spline model with unit error variance (Example~\ref{egpss})
and the shape-rate Gamma model with unit shape (Example~\ref{egsgm}),
which share the same joint limit distribution with an explicit
covariance matrix $\Psi$.

\begin{Corollary}[(Joint limit
distribution for partial smoothing spline model and shape-rate gamma model)]\label{l2regressionCLT}
Let $m>1+\sqrt{3}/2\approx1.866$,
and $h\asymp h^\ast$.
Suppose that (\ref{smoothG}) holds
for some $1\geq b>1/(2m)$, and $E(X-E(X|Z))^{\otimes2}$ is
positive definite. Furthermore, $g_0\in H^{m}(\mathbb{I})$
satisfies\break $\sum_\nu|V(g_0,h_\nu)|\nu^m<\infty$.
Then, as $n\rightarrow\infty$,
\begin{equation}
\label{psslimit}
\pmatrix{\sqrt{n}(\widehat{\theta}_{n,\lambda}-\theta_0)
\vspace*{3pt}\cr \sqrt{nh} \bigl(\widehat{g}_{n,\lambda}(z_0)-g_0(z_0)
\bigr)} \stackrel{d} {\longrightarrow} N (0,\Psi ),
\end{equation}
where
\[
\Psi=\pmatrix{
\bigl\{E \bigl[X-E(X|Z)
\bigr]^{\otimes2} \bigr\}^{-1} & 0
\vspace*{3pt}\cr
0& \displaystyle\frac{\int_0^\infty(1+x^{2m})^{-2}\,dx}{\pi}}.
\]
\end{Corollary}

In Corollary~\ref{l2regressionCLT}, we notice that the nonparametric
estimation bias asymptotically vanishes. This is due to the condition
$\sum_\nu|V(g_0,h_\nu)|\nu^m<\infty$, which imposes additional
smoothness on $g_0\in H^m(\mathbb I)$. Therefore, convergence rate
$n^{-m/(2m+1)}$ for $\widehat g_{n,\lambda}(z_0)$ is actually
sub-optimal given this additional smoothness (under $\lambda=\lambda
^\ast$). In practice, we select the smoothing parameter based on CV or
GCV; see~\cite{G02}.
%
\section{Joint hypothesis testing}\label{seclrt}

In this section, we propose likelihood ratio testing for a set of joint
local hypotheses in a general form (\ref{H0}). Under very general
conditions, the null limit distribution is proved to be a mixture of a
Chi-square distribution with $p$ degrees of freedom and a scaled
noncentral Chi-square distribution with one degree of freedom.
Obviously, these two Chi-square distributions are contributed by the
parametric and nonparametric components, respectively. Hence, we reveal
a new version of the Wilks phenomenon \cite{W38,FZZ01} which adapts to
the semi-nonparametric context. We further give more explicit null
limit distributions for three commonly used joint hypotheses. A key
technical tool used in this section is a \textit{restricted} version of JBR.

Consider the following joint hypothesis:
\begin{equation}
\label{H0}
H_{0}\dvtx  M\theta+Q g(z_0)=\alpha
\quad\mbox{vs.}\quad H_1\dvtx  M\theta+Q g(z_0)\neq\alpha,
\end{equation}
where $M=(M_1^T,\ldots,M_k^T)^T$ is a $k\times p$ matrix with $k\le
p+1$, $Q=(q_1,\ldots,q_k)^T$
and the $\alpha$ are $k$-vectors. Without loss of generality, we assume
$N\equiv(M, Q)$ to have elements in $\mathbb{I}=[0,1]$. We further
assume that the matrix $N$ has full rank. $M$, $Q$ and $\alpha$ are all
prespecified according to the testing needs. For example, when $N$ is
the identity matrix $I_{p+1}$ and $\alpha=(\theta_0^T, w_0)^T$, $H_0$
reduces to $(\theta^T, g(z_0))^T=(\theta_0^T, w_0)^T$. See
Corollary~\ref{remcase} for more examples. This provides another way
to construct the joint CIs for $(\theta_0^T, g_0(z_0))^T$ without
estimating $\Omega^{-1}$ or $\sigma_{z_0}$. The simultaneous testing of
two marginal hypotheses, that is, $H_0^P\dvtx  \theta=\theta_0$ and
$H_0^{N}\dvtx  g(z_0)=w_0$, can also be used for this purpose, but it
requires the very conservative Bonferroni correction. Moreover, our
joint hypothesis is more general, and the testing approach is more
straightforward to implement.

To define the likelihood ratio statistic, we establish the constrained
estimate under (\ref{H0}) in three steps:\vspace*{1pt} (i) arbitrarily choose
$(\theta^\dag, w^\dag)\in\mathbb{R}^p\times\mathbb{R}$
satisfying $M\theta^\dag
+Qw^\dag=\alpha$; (ii) define $\widehat{f}_{n,\lambda}^0\equiv
(\widehat
{\theta}_{n,\lambda}^0, \widehat{g}_{n,\lambda}^0)=\arg\max_{f\in
\mathcal{H}_0}
L_{n,\lambda}(f)$, where $\mathcal{H}_0\equiv\{(\theta,g)\in
\mathcal
{H}| M\theta+Qg(z_0)=0\}$ and
\begin{equation}\label{lnlam}
\quad L_{n,\lambda}(f)
%
= n^{-1}\sum_{i=1}^n\ell
\bigl(Y_i;X_i^T\theta+g(Z_i)+X_i^T
\theta^\dag +w^\dag \bigr) -(1/2)\lambda J(g,g);\hspace*{-10pt}
\end{equation}
(iii) define the constrained estimate as $\widehat{f}_{n,\lambda
}^{H_0}=(\widehat{\theta}_{n,\lambda}^0+\theta^\dag,
\widehat{g}_{n,\lambda}^0+w^\dag)$. Then, the LRT statistic is
$\mathrm{LRT}_{n,\lambda}=\ell_{n,\lambda}(\widehat{f}_{n,\lambda
}^{H_0})-\ell
_{n,\lambda}(\widehat{f}_{n,\lambda})$.

Given the inner product $\langle\cdot,\cdot\rangle$, we note that
$\mathcal{H}_0$ is a closed subset in $\mathcal{H}$ and thus a Hilbert
space. Hence, we will construct
the projections of the two operators $R_u$ and $P_\lambda$ (associated
with $\mathcal H$) onto the subspace $\mathcal{H}_0$, denoting them
$R_u^0$ and $P_\lambda^0$, respectively. Lemma~\ref{propextru} below
provides a preliminary step for the construction. Its proof is similar
to that of Proposition~\ref{proru} and is thus omitted.

\begin{Lemma}\label{propextru}
For any $u=(x,z)\in\mathcal U$ and $q\in\mathbb{I}$,
define
\[
H_{q,u}=(\Omega+\Sigma_\lambda)^{-1} \bigl(x-qA(z)
\bigr) \quad\mbox{and}\quad T_{q,u}=qK_z-A^TH_{q,u}.
\]
Let $R_{q,u}\equiv(H_{q,u}, T_{q,u})\in\mathcal{H}$. Then, for any
$f\in
\mathcal{H}$ and $u\in\mathcal{U}$, we have $\langle R_{q,u},\break  f\rangle
=x^T\theta+q g(z)$.
\end{Lemma}

Obviously, $R_{q,u}$ is a generalization of $R_u$ defined in
Proposition~\ref{proru}, that is, $R_u=R_{1,u}$. Lemma~\ref
{propextru} implies that the restricted parameter space $\mathcal
{H}_0$ can be rewritten as
\begin{equation}
\label{extru}
\mathcal{H}_0= \bigl\{f=(\theta,g)\in\mathcal{H}|
\langle R_{q_j,W_j},f\rangle =0, j=1,\ldots,k \bigr\},
\end{equation}
where $W_j=(M_j,z_0)$. Define
$H(Q,W)=(H_{q_1,W_1},\ldots,H_{q_k,W_k})$,
$T(Q,W)=(T_{q_1,W_1}(z_0),\ldots,T_{q_k,W_k}(z_0))$ and
$M_K=MH(Q,W)+QT(Q,W)$. Construct the projections
\[
R_u^0=R_u-\sum
_{j=1}^k\rho_{u,j}R_{q_j,W_j}
\quad\mbox{and} \quad P^0_\lambda f=P_\lambda f-\sum
_{j=1}^k\zeta_j(f)R_{q_j,W_j},
\]
where $(\rho_{u,1},\ldots,\rho_{u,k})^T=M_K^{-1}(MH_u+QT_u(z_0))$ and
$(\zeta_1(f),\ldots,\zeta_k(f))^T=M_K^{-1}(MH_{g}^*+QT_{g}^*(z_0))$.
Recall that $R_u\dvtx  u \mapsto(H_u, T_u)$ and $P_\lambda\dvtx  (\theta,g)
\mapsto(H_g^*,\break  T_g^*)$ in Proposition~\ref{proru}. The invertibility
of $M_K$ is given in the proof of Proposition~\ref{propKW} below.

Proposition~\ref{propKW} below says that $R_u^0$ and $P^0_\lambda$
defined above are indeed what we need.
\begin{Proposition}\label{propKW}
Let $f=(\theta,g)$ and $\tilde{f}=(\tilde{\theta},\tilde{g})$.
For any $u=(x,z)\in\mathbb{I}^p\times\mathbb{I}$,
$f,\tilde{f}\in\mathcal{H}_0$,
we have $\langle R_u^0, f\rangle=x^T\theta+g(z)$ and
$\langle P^0_\lambda f,\tilde{f}\rangle=\lambda J(g,\tilde{g})$.
\end{Proposition}

Based on Proposition~\ref{propKW}, we can write down the Fr\'{e}chet
derivatives of
$L_{n,\lambda}$ defined in (\ref{lnlam}) under $\mathcal{H}_0$ by
modifying those of $\ell_{n,\lambda}$ as follows:
replace $\theta$, $g$, $R_U$ and $P_\lambda$ by $\theta+\theta^\dag$,
$g+w^\dag$, $R_U^0$ and $P_\lambda^0$. For example,
\begin{eqnarray*}
&&DL_{n,\lambda}(f)\Delta f
\\
&&\qquad =\frac{1}{n}\sum_{i=1}^n \dot{
\ell}_a \bigl(Y_i;X_i^T
\theta+g(Z_i)+X_i^T\theta^\dag+w^\dag
\bigr) \bigl\langle R_{U_i}^0,\Delta f \bigr\rangle- \bigl
\langle P^0_\lambda f,\Delta f \bigr\rangle
\\
&&\qquad\equiv \bigl\langle S^0_{n}(f), \Delta f \bigr\rangle-
\bigl\langle P^0_\lambda f,\Delta f \bigr\rangle= \bigl\langle
S^0_{n,\lambda}(f), \Delta f \bigr\rangle.
\end{eqnarray*}
Similarly, we have $S^0_{n,\lambda}(\widehat{f}_{n,\lambda}^0)=0$. Also
define $S^0(f)=E\{S_n^0(f)\}$ and $S_\lambda^0(f)=S^0(f)-P^0_\lambda
(f)$. For the second derivative, we have
$\mathit{DS}_{n,\lambda}^0(f)\Delta f_1\Delta f_2= \break D^2L_{n,\lambda}(f)\Delta
f_1\Delta f_2$ and $\mathit{DS}_\lambda^0(f)\Delta f_1\Delta f_2=\mathit{DS}^0(f)\Delta
f_1\Delta f_2-\langle P^0_\lambda\Delta f_1,\Delta f_2\rangle$, where
\[
\mathit{DS}^0(f)\Delta f_1\Delta f_2=E \bigl\{\ddot{
\ell}_a \bigl(Y;X^T\theta +g(Z)+X^T\theta
^\dag+w^\dag \bigr) \bigl\langle R_U^0,
\Delta f_1 \bigr\rangle \bigl\langle R_U^0,
\Delta f_2 \bigr\rangle \bigr\}.
\]

In Theorem~\ref{LRTmainthm1} below, we present a new version of JBR
that is restricted to
the subspace $\mathcal H_0$. We need an additional Assumption~\ref{A.6}
here. Let
$f_0^0\equiv(\theta_0-\theta^\dag,g_0-w^\dag)$, which belongs to
$\mathcal{H}_0$ under $H_0$.
Assumption~\ref{A.6} holds under mild conditions similar to those
specified in Proposition~\ref{ratesconvergence}. The proof can be
similarly conducted by replacing
the space $\mathcal{H}$ by $\mathcal{H}_0$, and thus, is omitted.

\begin{Assumption}\label{A.6}
Under $H_0$ specified in (\ref{H0}),
\[
\bigl\|\widehat{f}_{n,\lambda}^0-f_0^0\bigr\|= O_P\bigl((nh)^{-1/2}+h^m\bigr).
\]
\end{Assumption}

\begin{Theorem}[(Restricted joint Bahadur
representation)]\label{LRTmainthm1}  Suppose that Assumptions \ref{A.1}, \ref{A.3},
\ref{A.2} and \ref{A.6} hold and that $h=o(1)$ and
$nh^2\rightarrow\infty$ as $n\rightarrow\infty$. Under $H_0$
specified in (\ref{H0}), we have
$\|\widehat{f}_{n,\lambda}^0-f_0^0-S_{n,\lambda}^0(f_0^0)\|
=O_{P}(a_n\log{n})$, where $a_n$ is defined as in
(\ref{bahareps}).
\end{Theorem}

Given the above preparatory results, we are ready to present general
results for the null limit distribution of $ -2n\cdot \mathrm{LRT}_{n,\lambda}$
in Theorem~\ref{LRTlimdist}. Define $r_n=(nh)^{-1/2}+h^m$, and let
\[
\Phi_\lambda = \Lambda N^TM_K^{-1}N
\Lambda^T,
\]
where
\[
\Lambda=\pmatrix{
(\Omega+\Sigma_\lambda)^{-1/2}&0
\vspace*{2pt}\cr
0&K(z_0,z_0)^{1/2}}
 \pmatrix{
I_p&-A(z_0)
\vspace*{2pt}\cr
0&1 }.
\]

\begin{Theorem}[(Joint local testing)]\label{LRTlimdist}
Suppose that Assumptions \ref{A.1}\break through~\ref{A.6} are satisfied,
there exists $b\in(1/(2m),1]$ such that $G_k$ satisfies
(\ref{smoothG}), and $h=o(1)$,
$nh^2\rightarrow\infty$, $n^{1/2}h^{m(1+b)}=o(1)$,
$r_n^2h^{-1/2}=o(a_n)$ and
$a_n=o(\min\{r_n, n^{-1}r_n^{-1}(\log{n})^{-1},
n^{-1/2}h^{1/2}(\log{n})^{-1}\})$, where $a_n$ is defined as in
(\ref{bahareps}). Furthermore, for any $z_0\in[0,1]$,
$\lim_{h\rightarrow0}\sqrt{n}(W_\lambda g_0)(z_0)/\sqrt
{K(z_0,z_0)}=c_{z_0}$,
$\lim_{h\rightarrow0}\Phi_\lambda=\Phi_0$, where $\Phi_0$ is a
fixed $(p+1)\times(p+1)$ positive semidefinite matrix, and
\begin{eqnarray}\label{LRTacondition}
\lim_{h\rightarrow0}h V(K_{z_0},K_{z_0}) &\rightarrow &
\sigma_{z_0}^2>0,
\\
\label{LRTbcondition}
\lim_{h\rightarrow0}E_Z \bigl\{B(Z)\bigl|K_{z_0}(Z)\bigr|^2
\bigr\}/K(z_0,z_0) &\equiv &  c_0
\in(0,1].
\end{eqnarray}
Under $H_0$ specified in (\ref{H0}), we obtain: \textup{(i)} $\|\widehat
{f}_{n,\lambda}-\widehat{f}_{n,\lambda}^{H_0}\|
=O_P(n^{-1/2})$; \textup{(ii)} $-2n\times \mathrm{LRT}_{n,\lambda}=n\|\widehat
{f}_{n,\lambda}-\widehat{f}_{n,\lambda}^{H_0}\|^2+o_{P}(1)$;
\begin{equation}
\mathrm{(iii)} -2n\cdot \mathrm{LRT}_{n,\lambda}\stackrel{d} {\longrightarrow}\upsilon
^T\Phi _0\upsilon,\label{mainlim}
\end{equation}
where $\upsilon\sim N
\bigl({0 \choose c_{z_0}},
{I_p  \ \  0\choose  0 \ \   c_0}\bigr)$.
\end{Theorem}

The \textit{parametric} convergence-rate result proved in (i) of
Theorem~\ref{LRTlimdist} is reasonable since the null hypothesis
imposes only a finite-dimensional constraint. By~(\ref{PropRPK}), it
can be explicitly
shown that
\begin{equation}
\label{findingc0}
\quad c_0=\lim_{\lambda\rightarrow0}\frac{Q_2(\lambda,z_0)}{
Q_1(\lambda,z_0)}\qquad
\mbox{where } Q_l(\lambda,z)\equiv\sum_{\nu\in
\mathbb{N}}
\frac{|h_\nu(z)|^2}{(1+\lambda\gamma_\nu)^l} \mbox{ for }l=1,2.\hspace*{-7pt}
\end{equation}
It is well known that the reproducing kernel $K$ is uniquely
determined for any Hilbert space if it exists; see \cite{S97}, page 38.
This implies that $c_0$ defined in (\ref{LRTbcondition}) is also
uniquely determined. Therefore, different choices of $(h_\nu, \gamma
_{\nu})$ in (\ref{findingc0}) will give exactly the same value of
$c_0$, although a particular choice may facilitate the computation of
$c_0$. For example, in case (I) of Example~\ref{egpss}, we can
explicitly calculate $c_0$ as $0.75$ (0.83) when $m=2$ (3) by
choosing the trigonometric basis~(\ref{H0basis}).

The null limit distribution derived in Theorem~\ref{LRTlimdist}
cannot be directly used for inference because of the nontrivial
estimation of $c_{z_0}$. Hence, in Corollary~\ref{LRTcor2}, we present
a set of conditions under which the estimation bias of $\widehat
g_{n,\lambda}$ can be removed, and thus $c_{z_0}=0$.

\begin{Corollary}\label{LRTcor2}
Suppose that Assumptions \ref{A.1} through \ref{A.6} are satisfied,
and hypothesis $H_0$ holds. Let $m>1+\sqrt{3}/2\approx1.866$ and
$G_1,\ldots,G_p$ satisfy (\ref{smoothG}) with $1/(2m)<b\le1$. Also
assume that the Fourier coefficients
$\{V(g_0,h_\nu)\}_{\nu\in\mathbb{N}}$ of $g_0$ satisfy
$\sum_\nu|V(g_0,h_\nu)| \gamma_\nu^{1/2}<\infty$.
Furthermore, if $\Phi_\lambda$ converges to some fixed $(p+1)\times
(p+1)$ positive semidefinite matrix, that is, $\Phi_0$, and (\ref
{LRTacondition}) and (\ref{LRTbcondition})
are both satisfied\vspace*{1pt} for any $z_0\in[0,1]$, then (\ref{mainlim}) holds
with $c_{z_0}=0$
given that $h=h^\ast\asymp n^{-1/(2m+1)}$.
\end{Corollary}

Combining Theorem~\ref{LRTlimdist} with Corollary~\ref{LRTcor2}, we
immediately obtain Corollary~\ref{remcase}, which gives null limit
distributions of the three commonly assumed joint hypotheses.

\begin{Corollary}\label{remcase}
Suppose that the conditions in Corollary~\ref{LRTcor2} hold. We have:
\begin{longlist}[(III)]
\item[(I)] $H_0\dvtx  \theta=\theta_0$ and $g(z_0)=w_0$:
\[
-2n\cdot\mathrm{LRT}_{n,\lambda}\stackrel{d} {\longrightarrow} \chi_p^2+c_0
\chi_1^2,
\]
where the two Chi-square distributions are independent. In this case,
$N=I_{p+1}$, $\alpha=(\theta_0^T,w_0)^T$ and $\Phi_\lambda=\Phi
_0=I_{p+1}$.

\item[(II)] $H_0\dvtx  D\theta=\theta_{0}'$ and
$g(z_0)=w_0$ [$D$ is an $r\times p$ matrix with $0<r\le p$ and $\operatorname{rank}(D)=r$,
$\theta_{0}'$
is an $r$-vector with $0<r<p$]:
\[
-2n\cdot \mathrm{LRT}_{n,\lambda}\stackrel{d} {\longrightarrow}\chi
_r^2+c_0\chi_1^2,
\]
where the two Chi-square distributions\vspace*{-2pt} are independent. In this case,
$N={D \ \ 0_r\choose  0_p^T \ \  1}$,
$\alpha=(\theta_{0}^{\prime T},w_0)^T$ and
$\Phi_0={ \mathcal{P}_r  \ \  0_p\choose
0_p^T \ \  1}$ with\vspace*{1.5pt} the projection matrix (of rank~$r$)
$\mathcal{P}_r=\Omega^{-1/2}D^T(D\Omega^{-1}D^T)^{-1}D\Omega^{-1/2}$.

\item[(III)]
$H_0\dvtx  x_0^T\theta+g(z_0)=\alpha$ ($\alpha$, $x_0$ and $z_0$ are given):
\[
-2n\cdot \mathrm{LRT}_{n,\lambda}\stackrel{d} {\longrightarrow} c_0
\chi_1^2.
\]
In this case, $N=(x_0^T,1)$ and
$\Phi_0={
0_{p\times p} \ \ 0_p\choose
0_p^T  \ \ \   1}$.
\end{longlist}
\end{Corollary}

The independence between the two Chi-square distributions in
(I) and (II) follows from the joint asymptotics phenomenon that
$\widehat\theta_{n,\lambda}$ and $\widehat g_{n,\lambda}(z_0)$ are
asymptotically independent. In comparison with (I) and (II), we note
that the null limit distribution in (III) is dominated by the effect
from $g(z_0)$ because of its nonparametric nature, that is, its slower
convergence rate.

As far as we are aware, Corollary~\ref{remcase} is a new version of
the Wilks phenomenon \cite{W38,FZZ01} that adapts to the
semi-nonparametric context. Note that the value of $c_0$ converges to
one as $m\rightarrow\infty$. Therefore, this new type of Wilks
phenomenon reverts to the classical version in the parametric setup as
$m\rightarrow\infty$ by further consideration of the independence of
the two Chi-squares. For example, the null limit distribution in (I) of
Corollary~\ref{remcase} becomes $\chi^2_{p+1}$ as $m\rightarrow
\infty$.

In the end of this section, we apply Theorem~\ref{LRTlimdist}
to partial smoothing spline models (Example~\ref{egpss})
and shape-rate gamma models (Example~\ref{egsgm}). For simplicity, let
$0<z_0<1$.
Corollary~\ref{l2pssgammaLRT} directly follows from Corollary~\ref{LRTcor2}
and equivalent kernel theory \cite{MG93,N95}.

\begin{Corollary}[(Joint local testing for partial smoothing spline model and
shape-rate gamma model)]\label{l2pssgammaLRT}
Suppose that the hypothesis $H_0$ specified in \textup{(I)}
[\textup{(II)} or \textup{(III)}]
in Corollary~\ref{remcase} holds.
Let $m>1+\sqrt{3}/2\approx1.866$,
$G_1,\ldots,G_p$ satisfy (\ref{smoothG}) with $1/(2m)<b\le1$ and
$h\asymp h^\ast$. Also
assume that $g_0\in H^m(\mathbb{I})$ satisfies
$\sum_\nu|V(g_0,h_\nu)| \nu^m<\infty$, and
$E(X-E(X|Z))^{\otimes2}$ is positive definite.
Then, as $n\rightarrow\infty$,
the conclusion of \textup{(I)}
[\textup{(II)} or \textup{(III)}]
in Corollary~\ref{remcase} holds
with $c_0=\frac{\pi(z_0)\int_\mathbb{R}\omega_0(t)^2\,dt}{\omega_0(0)}$,\vspace*{1.5pt}
where the equivalent kernel function $\omega_0$ is specified in \cite{MG93}, page 184.
In particular, when $m=2$ (3) and the design is uniform,
$c_0=0.75$ (0.83).
\end{Corollary}

As for the logistic regression model (Example~\ref{exlog}), we need to
numerically approximate the value of $c_0$ due to the implicit forms of
the eigenfunctions and eigenvalues; see more detailed discussions in
Section~S.15 of the supplementary file \cite{CS15}.

\section{Examples}\label{secexa}

In this section, we present three concrete examples together with
simulations. In all the examples, the $G_k$s are sufficiently smooth
for Theorem~\ref{mainthm3} and Corollary~\ref{remcase} to apply.
Detailed assumption verifications for three examples can be found in
Sections~S.9,~S.13 and~S.15 of \cite{CS15}.

\begin{Example}[(Partial smoothing spline)]\label{egpss}
Consider a partially linear regression model
\begin{equation}
\label{npmodel}
Y=X^T\theta+g(Z)+\epsilon,
\end{equation}
where $\epsilon\sim N(0,\sigma^2)$ with an unknown $\sigma^2$. Hence,
$B(Z)=\sigma^{-2}$. For simplicity, $Z$ is assumed to be uniformly
distributed over $\mathbb I$. In this case, $V(g,\tilde{g})$
becomes the usual $L^2$-norm. The function \textit{ssr}() in the R package
\textit{assist} was used to select the smoothing parameter $\lambda$
based on CV or GCV; see \cite{KW02}. The unknown error variance can be
consistently estimated by $\widehat{\sigma}^2=n^{-1}\sum_i
(Y_i-X_i^T\widehat{\theta}_{n,\lambda}-\widehat{g}_{n,\lambda
}(Z_i))^2/(n-\operatorname{trace}(A(\lambda)))$,
where $A(\lambda)$ denotes the smoothing matrix; see \cite{W90}.

We next consider two separate cases: (I) $g\in H_0^m(\mathbb{I})$ and
(II) $g\in H^m(\mathbb{I})$.

\begin{longlist}[{}]
\item[\textit{Case} (I) $g\in H_0^m(\mathbb{I})$:] We choose the
following trigonometric eigensystem for $H_0^m(\mathbb{I})$:\vspace*{-2pt}
\begin{equation}
\label{H0basis}
h_\mu(z)=\cases{
\sigma,& \quad$\mu=0,$
\vspace*{2pt}\cr
\sqrt{2}\sigma\cos(2\pi kz), &$\quad\mu=2k, k=1,2,\ldots,$
\vspace*{2pt}\cr
\sqrt{2}\sigma\sin(2\pi kz), &$\quad\mu=2k-1, k=1,2,\ldots,$}
\end{equation}
with the corresponding $\gamma_\nu$ specified as
$\gamma_0=0$ and $\gamma_{2k-1}=\gamma_{2k}=\sigma^2(2\pi k)^{2m}$
for $k\ge1$.

It follows from (\ref{sigma}) and (\ref{H0basis}) that the asymptotic
variance of $\widehat{g}_{n,\lambda}(z_0)$ is expressed\vspace*{-3pt} as
\[
\sigma_{z_0}^2=\lim_{h\rightarrow0} \Biggl\{
\sigma^2 h \Biggl(1+2 \sum_{k=1}^{\infty}
\frac{1}{(1+(2\pi
h \sigma^{1/m}k)^{2m})^2} \Biggr) \Biggr\}.
\]
Lemma~6.1 in \cite{SC12} leads to,\vspace*{-2pt} for $l=1,2$,
\begin{equation}
\label{pssaeq}
\sum_{k=1}^{\infty}
\frac{1}{(1+(2\pi
h \sigma^{1/m}k)^{2m})^l}\sim\frac{I_l}{2\pi h \sigma^{1/m}},
\end{equation}
where $I_l=\int_0^\infty(1+x^{2m})^{-l}\,dx$. Therefore, we have
$\sigma
_{z_0}^2=(I_2\sigma^{2-1/m})/\pi$.
According to Corollary~\ref{l2regressionCLT}, the 95\% prediction
interval for $Y$ at a new observed covariate
$u_0=(x_0,z_0)$\vspace*{-2pt} is
\begin{equation}
\label{predintformula}
\widehat{Y}\pm1.96 \sqrt{\widehat{\sigma}^{2-1/m}I_2/(
\pi nh)+\widehat {\sigma}^2},
\end{equation}
where $\widehat{Y}=x_0^T\widehat{\theta}_{n,\lambda}+\widehat
{g}_{n,\lambda}(z_0)$
is the predicted response. We next calculate $c_0$ based on (\ref
{findingc0}). It follows from (\ref{H0basis}) and (\ref{pssaeq}) that
\begin{eqnarray*}
Q_l(\lambda,z_0)&=&\sigma^2+\sum
_{k\ge1} \biggl\{\frac
{|h_{2k}(z_0)|^2}{(1+\lambda\sigma^2(2\pi
k)^{2m})^l}+\frac{|h_{2k-1}(z_0)|^2}{(1+\lambda\sigma^2(2\pi
k)^{2m})^l} \biggr\}
\\[-2pt]
&=&\sigma^2+2\sigma^2\sum_{k\ge1}
\frac{1}{(1+\lambda\sigma
^2(2\pi
k)^{2m})^l}
\\[-2pt]
&=&\sigma^2+2\sigma^2\sum_{k\ge1}
\frac{1}{(1+(2\pi h \sigma
^{1/m}k)^{2m})^l}\sim\frac{I_l}{\pi h \sigma^{1/m}}
\end{eqnarray*}
for $l=1,2$. Hence we\vspace*{-3pt} obtain
\begin{equation}
\label{explicitc0}
c_0=I_2/I_1.
\end{equation}
Further calculations reveal that $c_0=0.75$ (0.83) when $m=2$ (3).

In the simulations, we first verify the joint asymptotics phenomenon,
that is, (\ref{animportantlimit}),
by investigating\vspace*{1pt} the (asymptotic) independence between $\widehat
{\theta
}_{n,\lambda}$ and $\widehat{g}_{n,\lambda}(z_0)$. Let $\theta_0=(8,-8)^T$
and $g_0(z)=0.6\beta_{30,17}(z)+0.4\beta_{3,11}(z)$, where $\beta
_{a,b}$ is the density function for $\operatorname{Beta}(a,b)$. We estimate the
nonparametric function $g_0$, which has many peaks and troughs, using
periodic splines with $m=2$; $\sigma$ is set to one. To allow the
linear and nonlinear\vspace*{-1pt} covariates $(X,Z)$ to be dependent, we generate
them as follows: generate $U, V, Z\stackrel{\mathrm{i.i.d.}}{\sim}\operatorname{Unif}[0,1]$, and
set $X_1=(U+0.2Z)/1.2$, $X_2=(V+0.2Z)/1.2$.
This leads to $\operatorname{corr}(X_1,Z)=\operatorname{corr}(X_2,Z)\approx0.20$, where $\operatorname{corr}$
denotes the correlation coefficient.
The dependence between $\widehat{\theta}_{n,\lambda}$ and $\widehat
{g}_{n,\lambda}(z)$ is evaluated through the absolute values of the
sample correlation coefficients (ACC)
between $\widehat{\theta}_{n,\lambda}=(\widehat\theta_{n,\lambda,1},
\widehat\theta_{n,\lambda,2})^T$ and $\widehat{g}_{n,\lambda}(z)$ at
ten evenly spaced $z$ points in $[0,1]$ based on $500$ replicated data
sets. The results are summarized in Figure~\ref{corthetagcase1} for
sample sizes $n=100,300,1000$. As $n$ increases, it is easy to see that
the ACC curves become uniformly closer to zero, which strongly
indicates the desired asymptotic independence.
%
\begin{figure}[b]

\includegraphics{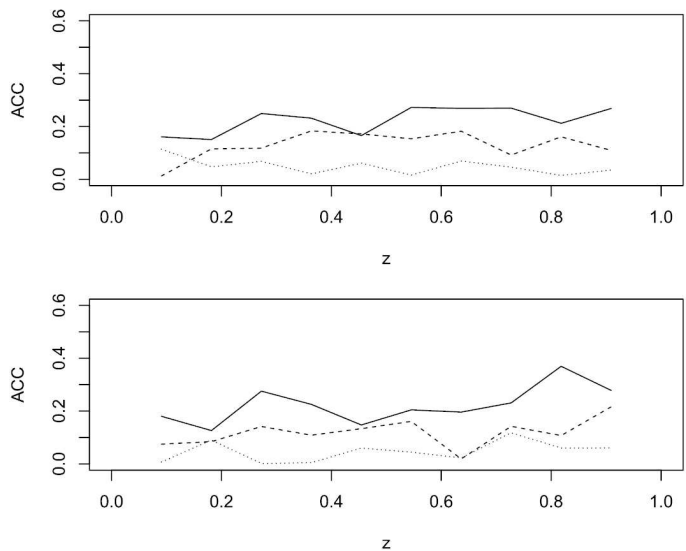}

\caption{Absolute values of correlation
coefficients (ACC)
between $\widehat{\theta}_{n,\lambda,1}$ and $\widehat
{g}_{n,\lambda
}(z)$ (the upper plot), and $\widehat{\theta}_{n,\lambda,2}$ and
$\widehat{g}_{n,\lambda}(z)$
(the lower plot), at ten evenly spaced nonlinear covariates in case \textup{(I)}
of Example~\protect\ref{egpss}.
The three lines correspond to three sample sizes: $n=100$ (solid),
$n=300$ (dashed), $n=1000$ (dotted).}
\label{corthetagcase1}
\end{figure}

To examine the performance of the $95\%$ prediction intervals (\ref
{predintformula}),
we calculate the proportions of the prediction intervals
covering the future response $Y$ generated from model (\ref{npmodel}),
that is, the coverage proportion. The simulation setup is the same as
before, except that we assume a one-dimensional $\theta_0=4$ for
simplicity. The new covariates are $(x_0,z_0)$ with $x_0=1/4,2/4,3/4$
and $z_0$ being thirty evenly spaced points in $[0,1]$. The coverage
proportions are calculated based on $500$ replications. We summarize
our simulation results in Figure~\ref{CPOfPIcase1} for sample sizes
$n=100,300,1000$. As $n$ grows, all the coverage proportions approach
the nominal level, $95\%$. In addition, the prediction interval lengths
approach the theoretical value indicated in formula (\ref
{predintformula}), that is, $2\times1.96=3.92$.
%
\begin{figure}

\includegraphics{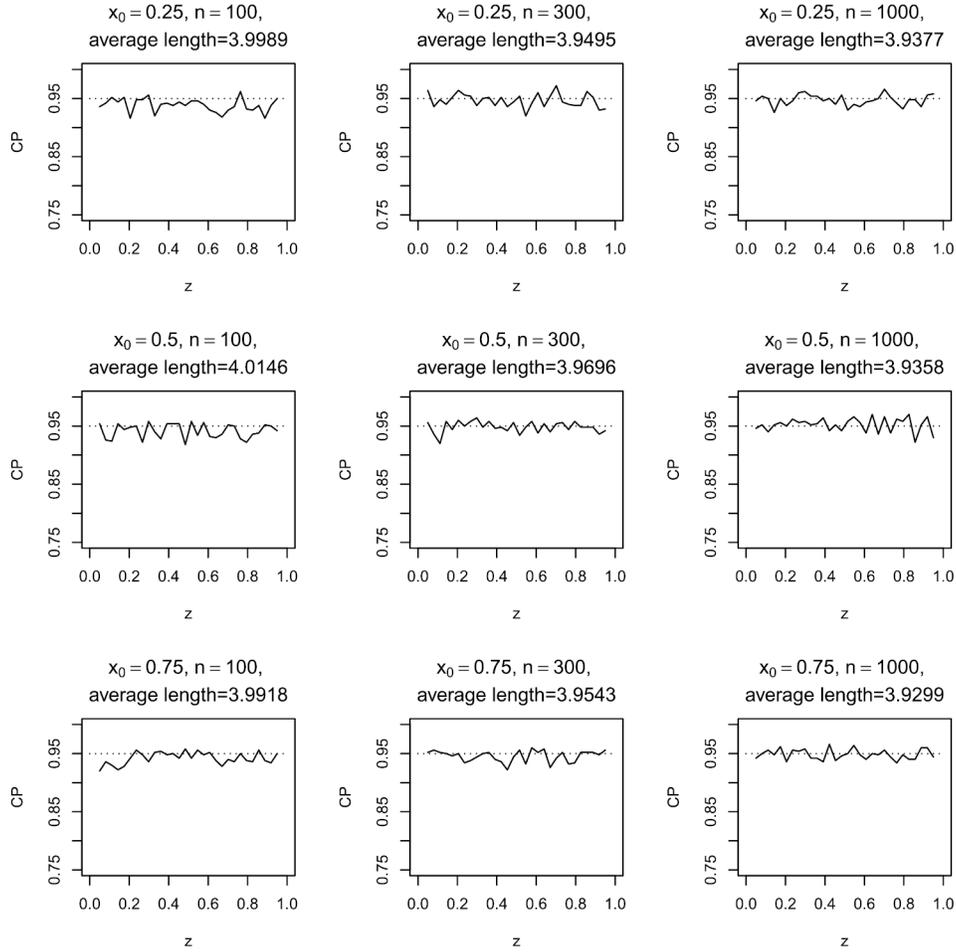}

\caption{Coverage proportion of $95\%$ prediction
intervals in case \textup{(I)} of Example~\protect\ref{egpss}.}
\label{CPOfPIcase1}
\end{figure}



Finally, we test $H_0\dvtx  x_0\theta+g(z_0)=0$. The true parameters are
chosen as $\theta_0=-4$, $g_0(z)=\sin(\pi z)$ and $\sigma=1$. The
performance is demonstrated by calculating the powers for the nine
combinations of $x_0=1/4,2/4,3/4$ and $z_0=1/4,2/4,3/4$ through $500$
replicated data sets.
In particular, $H_0$ is true when $x_0=1/4$ and $z_0=2/4$,
and $H_0$ is false at the other values of $(x_0,z_0)$. The results are
summarized in Table~\ref{powerloclrt} for sample sizes
$n=50,100,300,500,1000,1500$.
We observe that when $x_0=1/4$ and $z_0=2/4$, the power approaches the
correct size $5\%$,
while at the other values of $(x_0,z_0)$, where $H_0$ does not hold,
the power approaches one. This shows the validity of our local LRT test.
The detailed computational algorithm for the constrained estimate under
$H_0$ is given in Section~S.16 of the supplementary document \cite{CS15}.

\begin{table}
\caption{$100{\times}$ power of the local LRT test for nine combinations
of $x_0$ and $z_0$ for case \textup{(I)} of Example~\protect\ref{egpss}}
\label{powerloclrt}
\begin{tabular*}{\tablewidth}{@{\extracolsep{\fill
}}ld{2.2}d{3.2}d{3.2}d{3.2}d{3.2}d{3.2}@{}}
\hline
& \multicolumn{1}{c}{$\bolds{n=50}$} &\multicolumn{1}{c}{$\bolds
{n=100}$} &
\multicolumn{1}{c}{$\bolds{n=300}$} &
\multicolumn{1}{c}{$\bolds{n=500}$} &
\multicolumn{1}{c}{$\bolds{n=1000}$} &
\multicolumn{1}{c@{}}{$\bolds{n=1500}$}\\
\hline
$x_0=1/4$\\
\quad$z_0=1/4$ & 43.00 &56.60 &77.60 &90.40 &97.80
&98.60\\
\quad$z_0=2/4$ & 20.60 &13.00 &7.20 &7.00 &5.60 &5.10\\
\quad$z_0=3/4$ & 42.00 &50.00 &77.60 &89.60 &97.80 &99.20\\[3pt]
$x_0=2/4$ \\
\quad $z_0=1/4$ & 98.60 &99.80 &100 &100 &100 &100\\
\quad$z_0=2/4$ & 96.80 &99.00 &100 &100 &100 &100\\
\quad$z_0=3/4$ & 98.80 &99.80 &100 &100 &100 &100\\[3pt]
$x_0=3/4$\\
\quad$z_0=1/4$ & 99.80 &100 &100 &100 &100 &100\\
\quad$z_0=2/4$ & 99.60 &100 &100 &100 &100 &100\\
\quad$z_0=3/4$ & 99.60 &100 &100 &100 &100 &100\\
\hline
\end{tabular*}
\end{table}

\item[\textit{Case} (II) $g\in H^m(\mathbb{I})$:] For this larger
parameter space, we first construct an effective eigensystem that
satisfies (\ref{eigenproblem}). Let $\tilde{h}_\nu$s and $\tilde
{\gamma}_\nu$s be the normalized (with respect to the usual $L_2$-norm)
eigenfunctions and eigenvalues of the boundary value problem
$(-1)^m\tilde{h}_\nu^{(2m)}=\tilde{\gamma}_\nu\tilde{h}_\nu$,
$\tilde
{h}_\nu^{(j)}(0)=\tilde{h}_\nu^{(j)}(1)=0$, $j=m,m+1,\ldots,2m-1$. Then
we can construct $h_\nu=\sigma\tilde{h}_\nu$ and $\gamma_\nu
=\sigma^2
\tilde{\gamma}_\nu$.
Consequently,
\begin{eqnarray}
Q_l(\lambda,z)&=&\sum_{\nu}
\frac{|h_\nu(z)|^2}{(1+\lambda\gamma
_\nu
)^l}
\nonumber
\\[-8pt]
\label{npc0}
\\[-8pt]
\nonumber
&=&\sigma^{2-1/m}h^{-1}\sum_{\nu}
\frac{h\sigma^{1/m} |\tilde
{h}_\nu
(z)|^2}{(1+(h\sigma^{1/m})^{2m}\tilde{\gamma}_\nu)^l} \sim\sigma^{2-1/m}h^{-1} c_l(z),
\end{eqnarray}
where $c_l(z)=\lim_{h^\dag\rightarrow0}
\sum_{\nu}\frac{h^\dag|\tilde{h}_\nu(z)|^2}{(1+(h^{\dag
})^{2m}\tilde
{\gamma}_\nu)^l}$ and $h^\dag=h\sigma^{1/m}$, for $l=1,2$. Hence, by~(\ref{findingc0}), we have $c_0=c_2(z_0)/c_1(z_0)$. In addition, by
(\ref{sigma}), we obtain the asymptotic variance of $\widehat
{g}_{n,\lambda}(z_0)$ as $\sigma^{2-1/m} c_2(z_0)$, implying the
following 95\% prediction interval:
\[
\widehat{Y}\pm1.96 \sqrt{\widehat{\sigma }^{2-1/m}c_2(z_0)/(nh)+
\widehat {\sigma}^2}.
\]
The above discussion applies to general $m$. However, when $m=2$, we
can avoid estimating the $c_l(z_0)$s required in the inference by
applying the equivalent kernel approach. Following the discussion in
\cite{SC12}, we can actually obtain the same values of $c_0$ and
$\sigma
_{z_0}^2$ as in case (I). The simulation setup is the same as before
except that a different (nonperiodic) $g_0(z)=\sin(2.8\pi z)$ is used.
Figure~\ref{CPOfPIcase2} displays the coverage proportion of the 95\%
prediction intervals for three sample sizes $n=100,300,1000$. As $n$
grows, all the coverage proportions approach the $95\%$ nominal level,
and the prediction interval lengths approach the theoretical value~$3.92$.

\begin{figure}

\includegraphics{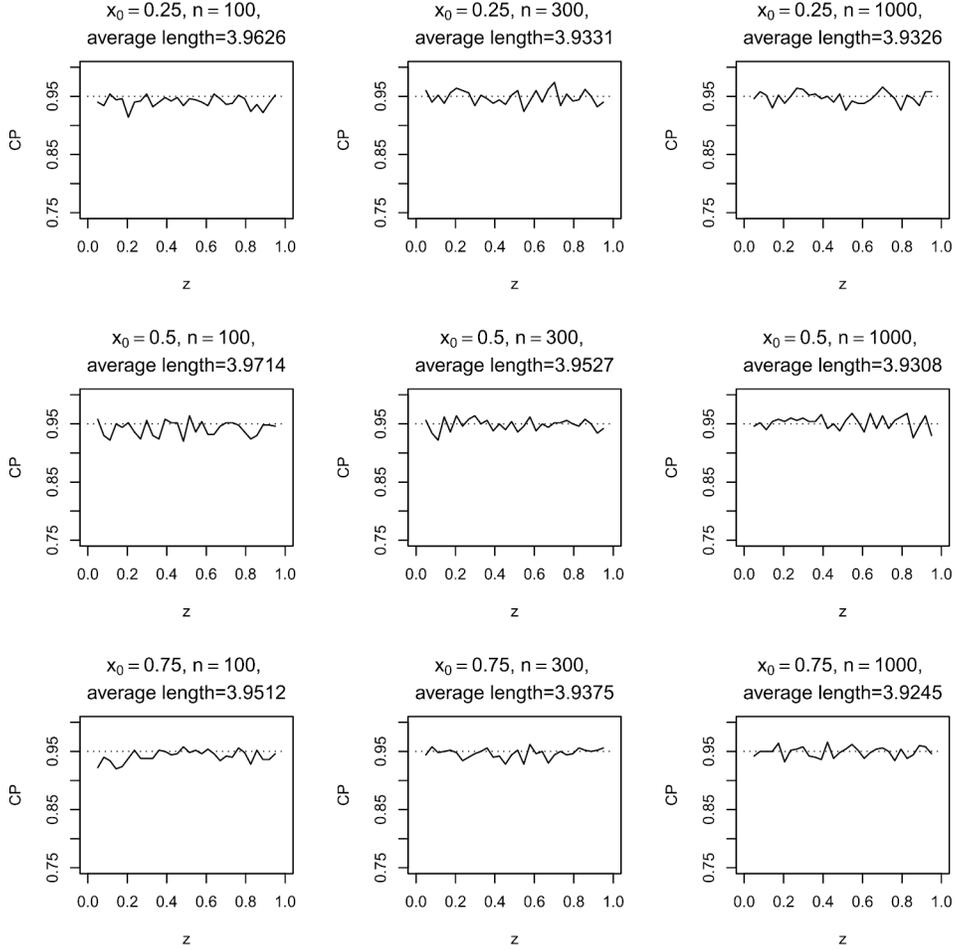}

\caption{Coverage proportion of 95\% prediction
intervals in case \textup{(II)} of Example~\protect\ref{egpss}.}
\label{CPOfPIcase2}
\end{figure}
%
\end{longlist}
\end{Example}

\begin{Example}[(Semiparametric gamma model)]\label{egsgm}
Consider a two-parameter exponential model
\[
Y| X,Z\sim\operatorname{Gamma} \bigl(\alpha,\exp \bigl(X^T
\theta_0+g_0(Z) \bigr) \bigr),
\]
where $\alpha>0$ is known, $g_0\in H_0^m(\mathbb{I})$ and $Z\sim
\operatorname{Unif}[0,1]$. It can be easily shown that $I(U)=\alpha$, and thus
$B(Z)=\alpha$ in this model. Consequently, we can construct the basis
functions $h_\nu$ as those defined in (\ref{H0basis}) with $\sigma
=\alpha^{-1/2}$, and the eigenvalues as $\gamma_0=0$ and $\gamma
_{2k-1}=\gamma_{2k}=\alpha^{-1}(2\pi k)^{2m}$ for $k\ge1$. The
remaining analysis is similar to case (I) of Example~\ref{egpss}; for
example, $c_0$ is given in (\ref{explicitc0}).
\end{Example}

\begin{Example}[(Semiparametric logistic regression)]\label{exlog}
For the binary response $Y\in\{0,1\}$, we consider the
following semiparametric logistic model:
\begin{equation}
\label{lgm}
P(Y=1|X=x,Z=z)=\frac{\exp(x^T\theta_0+g_0(z))}{1+\exp(x^T\theta_0+g_0(z))},
\end{equation}
where $g_0\in H^m(\mathbb I)$. It can be shown that, in reasonable situations,
all the conditions in Theorems \ref{mainthm3} and \ref{LRTlimdist}
are satisfied;
see Section~S.15 in \cite{CS15} for more details.

The solutions $\gamma_\nu$ and $h_\nu$ to the problem
(\ref{eigenproblem}) are useful to calculate the
quantities in the limit distribution (such as $\sigma_{z_0}^2$ and $c_0$
in Theorems \ref{mainthm2} and \ref{LRTlimdist}).
However, in this model, due to the intractable forms of these solutions,
we need to use consistent estimators of
$B(\cdot)$ and $\pi(\cdot)$ to find the approximated solutions; for
example, $\widehat B(\cdot)$ is a plug-in estimator and
$\widehat\pi(\cdot)$ is a kernel density estimator.

Given the length of this paper, we conduct simulations only for the CIs
of the conditional mean
defined in (\ref{lgm}) at a number of $(x_0,z_0)$ values, that is,
$x_0=1/4,2/4,3/4$ and thirty evenly
spaced $z_0$ over $[0,1]$.
The true parameters are $\theta_0=-0.5$ and
$g_0(z)=0.3(10^6)(1-z)^6+(10^4)(1-z)^{10}-2$.
For simplicity, we generate $X,Z\stackrel{\mathrm{i.i.d.}}{\sim}\operatorname{Unif}[0,1]$. Based
on $500$ replicated data
sets, we construct the 95\% CIs and calculate their coverage
proportions. The results are
summarized in Figure~\ref{logitCP} for various sample sizes
$n=400,500,700$. We observe that, as $n$
increases, the coverage proportions
approach the desired level, 95\%, and the lengths of the CI approach zero.
\begin{figure}

\includegraphics{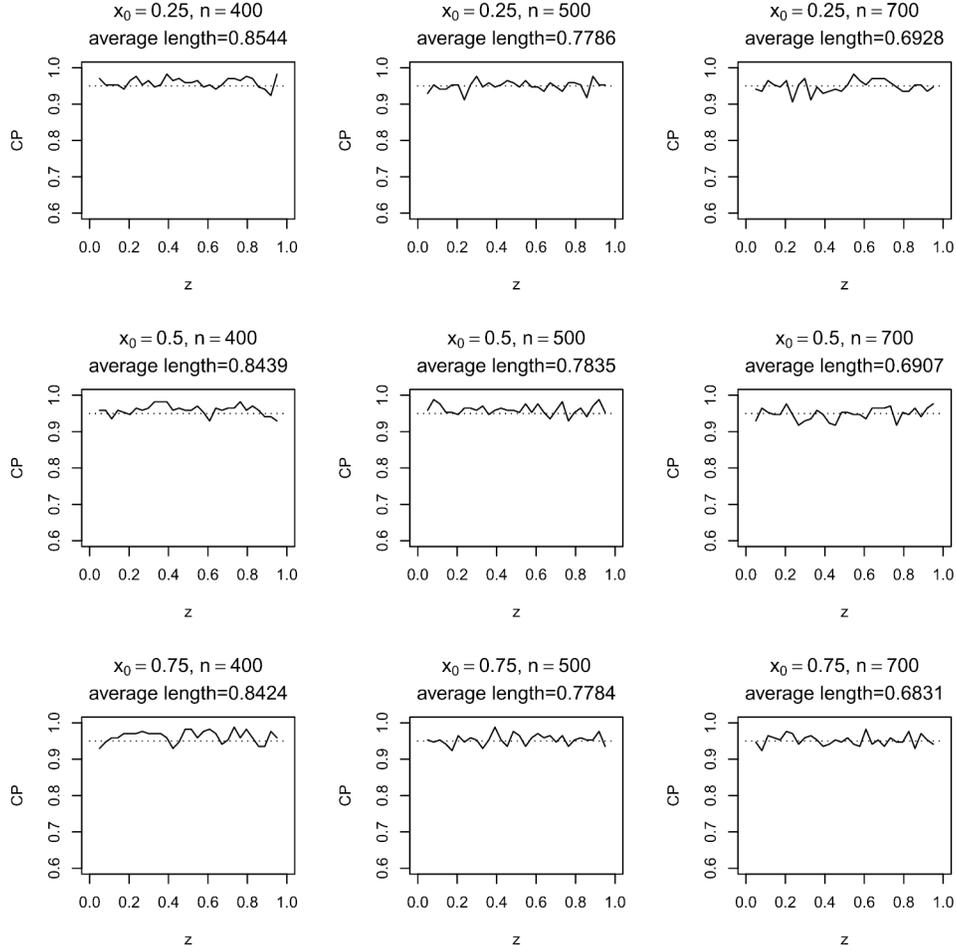}

\caption{Coverage proportion of 95\% CIs for the
conditional mean
constructed at a variety of $(x,z)$ values.}
\label{logitCP}
\end{figure}
\end{Example}

\begin{Remark}\label{remeg}
We use this logistic regression model to illustrate the eigensystem
difference between the semi-nonparametric context and the nonparametric
context, which leads to different inference for the nonparametric
components [except under some strong conditions, e.g., (\ref{equiva})
below]. This is slightly counterintuitive given that the parametric
component can be estimated at a faster rate. As discussed above, the
eigensystem for the semiparametric logistic model relies on $B(z)$
defined in (S.18) of \cite{CS15}. According to Shang and Cheng \cite
{SC12}, the eigensystem for the nonparametric\vspace*{1pt} logistic model relies on
$I'(z)$ defined as $\exp(g_0(z))/(1+\exp(g_0(z)))^{2}$. Therefore, the
equivalence of the two eigensystems holds if and only if $B(z)=I'(z)$,
that is,
\begin{equation}\label{equiva}
E \biggl\{\frac{\exp(X^T\theta_0)}{(1+\exp(X^T\theta
_0+g_0(z)))^2} \Big|Z=z \biggr\}= \frac{1}{(1+\exp(g_0(z)))^2}.
\end{equation}
If $\theta_0=0$, it is clear that (\ref{equiva}) is true. However, we
argue that in general (\ref{equiva}) may not hold. For instance, it
does not hold when $g_0(z)=0$ for\vspace*{1pt} some $z\in[0,1]$ because the above
equation then simplifies to
$E \{\frac{\exp(X^T\theta_0)}{(1+\exp(X^T\theta_0))^2}
\}=1$.
This is not possible since $\frac{\exp(X^T\theta_0)}{(1+\exp
(X^T\theta
_0))^2}<1$ almost surely.
This concludes our argument. \qed
\end{Remark}

\section{Future work}\label{secfut}

The general framework in this paper covers a wide range of commonly
used models. In this section, we discuss some possible extensions using
heuristic arguments, while omitting all the technical details due to
the length of this paper. The first possible extension is to the class
of generalized additive partial linear models in which
$E(Y|X,Z)=F(X^T\theta_0+\sum_{j=1}^{J}g_{j0}(Z_j))$. Our techniques are
expected to handle this more general class by modifying the joint
Bahadur representation in Theorem~\ref{mainthm1}, that is, to replace
$f$ therein by $(\theta, g_1,\ldots,g_J)$. The second possible
extension is to deal with the functional data. In \cite{SC142}, we
develop nonparametric inference for the (generalized) functional linear
models, that is, $E(Y|Z)=F(\int_0^1 Z(t)\beta_0(t)\,dt)$, by penalizing
the slope function $\beta(\cdot)$. By incorporating the techniques in
\cite{SC142} into our paper, we believe that it is feasible to do the
joint asymptotic study of the (generalized) partial functional linear
regression models \cite{S09}, that is, $E(Y| X,Z)=F(X^T\theta_0+\int_0^1 Z(t)\beta_0(t)\,dt)$. The third possible extension is from the above
regression models to survival models. Specifically, our results may be
extended to the partially linear Cox proportional hazard models (under
right censored data) (i.e., \cite{H99}), by replacing our criterion
function by their partial likelihood. This extension seems technically
feasible given the quadratic structure of the profile likelihood (a
generalization of partial likelihood) proven in~\cite{MV00}.

%

\begin{appendix}
\section*{Appendix}\label{app}

In this section, proofs of the main results are provided.
In Section~\ref{secanimportantlemma}, a preliminary lemma used for
main results
is provided. In Section~\ref{secproofthm2}, an initial result
about the joint limit distribution of the parametric and nonparametric
estimators with biased center
is given. Section~\ref{prooflemma4} includes the proof
of Theorem~\ref{mainthm3}. In Section~\ref{secproofjointlocaltest}, the
proof of Theorem~\ref{LRTlimdist} on the null limit distribution of
likelihood ratio testing is provided.

For any $f=(\theta,g)\in\mathcal{H}$, we treat $f$ as a ``partly linear''
function, that is, $f\dvtx (x,z)\mapsto x^T\theta+g(z)$, where $(x,z)\in
\mathcal{U}$.
Thus $(\theta,g)$ can be viewed as a bivariate function
defined on $\mathcal{U}$. Throughout the \hyperref[app]{Appendix}, we will not
distinguish $(\theta,g)$ and its associated function
$f$. For instance, we use
$(\theta,g)\in\mathcal{G}_0$ to mean
$f\in\mathcal{G}_0$, some set of functions defined over $\mathcal{U}$.

\subsection{An important lemma}\label{secanimportantlemma}
\begin{lemma}\label{lemma3}
\begin{eqnarray}\label{lemma3i}
\lim_{\lambda\rightarrow0}E_Z \bigl\{B(Z) \bigl(G(Z)-A(Z) \bigr)
\bigl(G(Z)-A(Z) \bigr)^T \bigr\}&=&0.
\\
\label{lemma3ii}
\lim_{\lambda\rightarrow0}E_Z \bigl\{B(Z)G(Z) \bigl(G(Z)-A(Z)
\bigr)^T \bigr\}&=& 0.
\end{eqnarray}
\end{lemma}

\begin{pf}
The proofs of (\ref{lemma3i}) and (\ref{lemma3ii}) are similar,
so we only show that (\ref{lemma3ii}) holds. Considering (\ref
{inter1}) and taking $g=h_\nu$, one has
\begin{equation}
\label{eqAG}
\quad V(G_k,h_\nu)=\langle
A_k,h_\nu \rangle_1= \biggl\langle \sum
_{\mu}V(A_k, h_\mu)h_{\mu},h_\nu
\biggr\rangle_1=(1+ \lambda\gamma _\nu
)V(A_k,h_\nu),\hspace*{-23pt}
\end{equation}
and, taking $g=K_z$, one has $V(G_k,K_z)=A_k(z)$. By (\ref{eqAG}),
$A_k=\sum_{\nu}\frac{V(G_k,h_\nu)}{1+\lambda\gamma_\nu}h_\nu$
holds in $L_{2}(\mathbb{I})$. For any $k,j=1,\ldots,p$, by a
straightforward calculation, we have
\[
E_Z \bigl\{B(Z)G_j(Z) \bigl(G_k(Z)-A_k(Z)
\bigr) \bigr\}=\sum_{\nu}V(G_j,h_\nu)V(G_k,h_\nu)
\frac{\lambda\gamma_\nu}{1+\lambda\gamma_\nu}.
\]
By square summability of $\{V(G_k,h_\nu)\}_{\nu\in\mathbb{N}}$ and
dominated convergence theorem, the above sum converges to zero as
$\lambda\rightarrow0$.
\end{pf}

\subsection{An initial result on joint asymptotic distribution with
biased center}\label{secproofthm2}

\begin{theorem}\label{mainthm2}
Let Assumptions \ref{A.1} through \ref{A.5} be satisfied.
Suppose that as $n\rightarrow\infty$, $h=o(1)$,
$nh^2\rightarrow\infty$ and $a_n\log{n}=o(n^{-1/2}h^{1/2})$,
where $a_n$ is defined as in (\ref{bahareps}). Furthermore,
assume that, as $n\rightarrow\infty$,
\begin{eqnarray}
 h V(K_{z_0},K_{z_0}) &\rightarrow &
\sigma^2_{z_0},\qquad  h^{1/2}(W_\lambda A)
(z_0)\rightarrow\alpha_{z_0}\in \mathbb{R}^p
\quad\mbox{and}
\nonumber
\\[-8pt]
\label{mainthm2part2}
\\[-8pt]
\nonumber
 h^{1/2}A(z_0) &\rightarrow &  -\beta_{z_0}
\in\mathbb{R}^p,
\end{eqnarray}
where $A$ is the Riesz representer defined in (\ref{inter1}). Then we
have, for any $z_0\in\mathbb{I}$,
\begin{equation}
\label{mainthm2eq} \pmatrix{
\sqrt{n} \bigl(\widehat{
\theta}_{n,\lambda}-\theta_0^* \bigr)
\vspace*{2pt}\cr
\sqrt{nh} \bigl(\widehat{g}_{n,\lambda}(z_0)-g_0^*(z_0)
\bigr) }
  \stackrel{d} {\longrightarrow} N
\bigl(0,\Psi^\ast \bigr),
\end{equation}
where
\begin{equation}
\label{Psieq}
\Psi^\ast=\pmatrix{
\Omega^{-1} & \Omega^{-1}(\alpha_{z_0}+\beta
_{z_0})
\vspace*{3pt}\cr
(\alpha_{z_0}+\beta_{z_0})^T\Omega^{-1}&
\sigma_{z_0}^2 +2\beta_{z_0}^T
\Omega^{-1}\alpha_{z_0}+\beta_{z_0}^T
\Omega ^{-1}\beta_{z_0} }.
\end{equation}
\end{theorem}

Note that $\Omega^{-1}$ is well defined under
Assumption~\ref{A.2}. It follows from (\ref{PropRPK}) and~(\ref{wlamA}) that
\begin{eqnarray*}
\alpha_{z_0} &=& \lim_{h\rightarrow0}h^{1/2}\sum
_{\nu}\frac{
V(G,h_\nu)\lambda\gamma_\nu}{(1+\lambda\gamma_\nu)^2}h_\nu
(z_0),
\\
\beta_{z_0} &=& -\lim_{h\rightarrow0}h^{1/2}\sum
_{\nu}\frac{V(G,h_\nu)}{1+\lambda\gamma_\nu}h_\nu(z_0).
\end{eqnarray*}

\begin{pf*}{Proof of Theorem~\ref{mainthm2}}
Define
\[
\widehat{f}_{n,\lambda}^h= \bigl(\widehat{\theta}_{n,\lambda
},h^{1/2}
\widehat {g}_{n,\lambda} \bigr),\qquad f_0^{*h}= \bigl(
\theta_0^*,h^{1/2}g_0^* \bigr),\qquad
R_u^h= \bigl(H_u,h^{1/2}T_u
\bigr),
\]
where we recall $f_0^*=(id-P_\lambda)f_0$, $H_u,T_u$ were defined by
(\ref{RTu}), and $P_\lambda$ is specified in Proposition~\ref{Penopr}.
By Theorem~\ref{mainthm1},
\[
\mathrm{Rem}_{n}=\widehat{f}_{n,\lambda}-f_0^*-
\frac{1}{n}\sum_{i=1}^n
\epsilon_i R_{U_i}
\]
satisfies $\|\mathrm{Rem}_n\|=O_{P}(a_n\log{n})$, which will imply by
Assumption~\ref{A.1}(b) that
\begin{equation}
\label{proofthm2eq1}
\Biggl\|\widehat{\theta}_{n,\lambda}-\theta_0^*-
\frac{1}{n}\sum_{i=1}^n
\epsilon_i H_{U_i}\Biggr\|_{l_2} =O_{P}(a_n
\log{n}).
\end{equation}
Define
$\mathrm{Rem}_n^h=\widehat{f}_{n,\lambda}^h-f_0^{*h}-\frac{1}{n}\sum_{i=1}^n
\epsilon_i R_{U_i}^h$, then it is easy to see that
\[
\Rem_n^h-h^{1/2} \Rem_n= \Biggl(
\bigl(1-h^{1/2} \bigr) \Biggl(\widehat{\theta}_{n,\lambda}-
\theta_0^*-\frac
{1}{n}\sum_{i=1}^n
\epsilon_i H_{U_i} \Biggr),0 \Biggr).
\]
Thus, by (\ref{proofthm2eq1}),
\begin{eqnarray*}
\bigl\|\Rem_n^h-h^{1/2} \Rem_n\bigr\|&\le&
\bigl(1-h^{1/2} \bigr)\cdot O \Biggl(\Biggl\|\widehat{\theta}_{n,\lambda}-
\theta_0^*-\frac{1}{n}\sum_{i=1}^n
\epsilon_i H_{U_i}\Biggr\|_{l_2} \Biggr)
\\
&=&O_{P}(a_n\log{n}).
\end{eqnarray*}
Since by assumption $a_n\log{n}=o(n^{-1/2})$,
$\|\Rem_n^h\|=o_{P}(n^{-1/2})$. Next we will use $\Rem_n^h$ to
obtain the target joint limiting distribution.

The idea is to employ the Cram\'{e}r--Wald device. For any $x\in
\mathbb{I}^p$,
we will obtain the limiting distribution of
$n^{1/2}x^T(\widehat{\theta}_{n,\lambda}-\theta
_0^*)+(nh)^{1/2}(\widehat
{g}_{n,\lambda}(z_0)-g_0^*(z_0))$.
Note that this is equal to $n^{1/2}\langle
R_u,\widehat{f}_{n,\lambda}^h-f_0^{*h}\rangle$ with $u=(x,z_0)$.
Using the fact that
\begin{eqnarray*}
&&\Biggl|n^{1/2} \Biggl\langle R_u, \widehat{f}_{n,\lambda}^h-f_0^{*h}-
\frac
{1}{n}\sum_{i=1}^n
\epsilon_i R_{U_i}^h \Biggr\rangle\Biggr|
\\
&&\qquad\le n^{1/2}\|R_u\|\cdot\bigl\|\Rem_n^h
\bigr\|
\\
&&\qquad=O_{P} \bigl(n^{1/2}h^{-1/2}a_n\log{n}
\bigr)=o_{P}(1),
\end{eqnarray*}
we just need to find the limiting distribution of $n^{1/2}\langle
R_u, \frac{1}{n}\sum_{i=1}^n \epsilon_i R_{U_i}^h\rangle$,
which is equal to
\[
n^{1/2} \Biggl\langle R_u, \frac{1}{n}\sum
_{i=1}^n \epsilon_i
R_{U_i}^h \Biggr\rangle=n^{-1/2}\sum
_{i=1}^n\epsilon_i
\bigl(x^T H_{U_i}+h^{1/2} T_{U_i}(z_0)
\bigr).
\]
Next we will use CLT to find its limiting distribution. By
Assumption~\ref{A.1}(c), that is,
$E\{\epsilon^2|U\}= I(U)$, we have that
\begin{eqnarray*}
s_n^2&\equiv&\Var \Biggl(\sum
_{i=1}^n\epsilon_i
\bigl(x^T H_{U_i}+h^{1/2} T_{U_i}(z_0)
\bigr) \Biggr)
\\
&=&nE \bigl\{\epsilon^2 \bigl|x^T H_U+h^{1/2}
T_U(z_0)\bigr|^2 \bigr\}
\\
&=&nE \bigl\{E \bigl\{\epsilon^2|U \bigr\} \bigl|x^T
H_U+h^{1/2} T_U(z_0)\bigr|^2
\bigr\}
\\
&=&n E \bigl\{I(U)\bigl|x^T H_U+h^{1/2}
T_U(z_0)\bigr|^2 \bigr\}.
\end{eqnarray*}
A direct examination from (\ref{RTu}) shows that
\begin{eqnarray}
&&x^T H_U+h^{1/2}
T_U(z)
\nonumber
\\
&&\qquad = x^T(\Omega+\Sigma_\lambda)^{-1} \bigl(X-A(Z)
\bigr)+h^{1/2}K_Z(z_0)
\nonumber
\\[-8pt]
\label{lindeq1}\\[-8pt]
\nonumber
&&\qquad\quad{}-h^{1/2}
A(z_0)^T (\Omega+\Sigma_\lambda)^{-1}
\bigl(X-A(Z) \bigr)
\\
&&\qquad= h^{1/2}K_Z(z_0)+ \bigl(x-h^{1/2}
A(z_0) \bigr)^T(\Omega+\Sigma_\lambda)^{-1}
\bigl(X-A(Z) \bigr).
\nonumber
\end{eqnarray}

It follows by the proof of Lemma~\ref{lemma0} that $|K_Z(z_0)|=O(h^{-1})$.
On the other hand, for any $z\in\mathbb{I}$ and\vspace*{-2pt} $j=1,\ldots,p$,
\begin{eqnarray*}
\bigl|A_k(z)\bigr|&=&\Biggl|\sum_{\nu=1}^\infty
\frac{V(G_k,h_\nu)h_\nu}{1+\lambda
\gamma
_\nu}\Biggr|
\\[-2pt]
&\le& \biggl(\sum_{\nu}\bigl|V(G_k,h_\nu)\bigr|^2h_\nu(z)^2
\biggr)^{1/2} \biggl(\sum_\nu
\frac
{1}{(1+\lambda\gamma_\nu)^2} \biggr)^{1/2} \le C'_k
h^{-1/2},
\end{eqnarray*}
where $C'_k$ is free of $z$. Thus, by (\ref{lindeq1}),
there exists a constant $c'$ s.t.
$|x^T H_U+h^{1/2} T_U(z)|\le c'h^{-1/2}$, a.s.

Thus\vspace*{-3pt}
\begin{eqnarray}
&&\quad E \bigl\{I(U)\bigl|x^T H_U+h^{1/2}
T_U(z_0)\bigr|^2 \bigr\}
\nonumber
\\[-1pt]
&&\quad\qquad=hE \bigl\{I(U) \bigl|K_Z(z_0)\bigr|^2 \bigr\}
\nonumber
\\[-9pt]
\label{eqvar}
\\[-9pt]
\nonumber
&&\qquad\qquad{}+2 h^{1/2} \bigl(x-h^{1/2}A(z_0)
\bigr)^T (\Omega+\Sigma_\lambda)^{-1} E \bigl
\{I(U)K_Z(z_0) \bigl(X-A(Z) \bigr) \bigr\}\hspace*{-12pt}
\nonumber
\\[-1pt]
&&\qquad\qquad{}+ \bigl(x-h^{1/2}A(z_0) \bigr)^T E \bigl
\{I(U)H_UH_U^T \bigr\} \bigl(x-h^{1/2}A(z_0)
\bigr).\nonumber
\end{eqnarray}
Lemma~\ref{lemma3} tells us, as $\lambda\rightarrow0$,
$\Sigma_\lambda=E_Z\{B(Z)G(Z)(G(Z)-A(Z))^T\}\rightarrow0$.
It can be verified that
\begin{eqnarray*}
\label{tmpeq2}
&& E_U \bigl\{I(U)H_U
H_U^T \bigr\}
\nonumber
\\[-1pt]
&&\qquad =(\Omega+\Sigma_\lambda)^{-1}E \bigl\{I(U) \bigl(X-A(Z)
\bigr) \bigl(X-A(Z) \bigr)^T \bigr\}(\Omega +\Sigma
_\lambda)^{-1}
\nonumber
\\[-1pt]
&&\qquad =(\Omega+\Sigma_\lambda)^{-1} E \bigl\{I(U)
\bigl(X-G(Z)+G(Z)-A(Z) \bigr)
\nonumber
\\[-1pt]
&&\hspace*{60pt}\qquad\qquad{}\times \bigl(X-G(Z)+G(Z)-A(Z) \bigr)^T \bigr\}(\Omega+
\Sigma_\lambda)^{-1}
\nonumber
\\[-1pt]
&&\qquad =(\Omega+\Sigma_\lambda)^{-1} \bigl(E \bigl\{I(U)
\bigl(X-G(Z) \bigr) \bigl(X-G(Z) \bigr)^T \bigr\}
\nonumber
\\[-1pt]
&&\hspace*{57pt}\qquad\quad{}+E \bigl\{I(U) \bigl(G(Z)-A(Z) \bigr) \bigl(G(Z)-A(Z) \bigr)^T
\bigr\} \bigr) (\Omega+\Sigma_\lambda)^{-1}\\
&&\qquad\rightarrow
\Omega^{-1},
\end{eqnarray*}
where the last limit follows by Lemma~\ref{lemma3}. By assumption,
as $\lambda\rightarrow0$, $hE\{I(U) |K_Z(z_0)|^2\}=h
V(K_{z_0},K_{z_0})\rightarrow\sigma_{z_0}^2$,
$h^{1/2}A(z_0)\rightarrow-\beta_{z_0}$ and
\begin{eqnarray*}
&&h^{1/2}E \bigl\{I(U)K_Z(z_0) \bigl(X-A(Z)
\bigr) \bigr\}
\\[-1pt]
&&\qquad =h^{1/2}E \bigl\{B(Z)K_{z_0}(Z) \bigl(G(Z)-A(Z) \bigr) \bigr
\}
\\[-1pt]
&& \qquad =h^{1/2} \bigl(V(G,K_{z_0})-V(A,K_{z_0}) \bigr)
\\[-1pt]
&&\qquad =h^{1/2} \bigl(A(z_0)-V(A,K_{z_0}) \bigr)
\\[-1pt]
&&\qquad =h^{1/2} (W_\lambda A) (z_0)\rightarrow
\alpha_{z_0}.
\end{eqnarray*}
Thus, as $\lambda$ approaches zero, the limit of (\ref{eqvar}) is
\[
\sigma_{z_0}^2+2(x+\beta_{z_0})^T
\Omega^{-1}\alpha_{z_0}+(x+\beta _{z_0})^T
\Omega^{-1}(x+\beta_{z_0})= \bigl(x^T,1 \bigr)
\Psi^\ast \bigl(x^T,1 \bigr)^T,
\]
where $\Psi^\ast$ is defined in (\ref{Psieq}).
So $s_n^2\asymp n$. Then it can be shown that, for any $\varepsilon>0$,
\begin{eqnarray*}
&&E \bigl\{\bigl|\epsilon \bigl(x^T H_U+h^{1/2}
T_U(z_0) \bigr)\bigr|^2 I{\bigl|\epsilon
\bigl(x^T H_U+h^{1/2} T_U(z_0)
\bigr)\bigr|\ge\varepsilon s_n} \bigr\}
\\
&&\qquad \le \bigl(c'h^{-1/2} \bigr)^2E \bigl\{
\epsilon^2I \bigl(|\epsilon|\ge\varepsilon s_n
h^{1/2}/c' \bigr) \bigr\}
\\
&&\qquad \le \bigl(c'h^{-1/2} \bigr)^2E \bigl\{
\epsilon^4 \bigr\}^{1/2}P \bigl(|\varepsilon|\ge \varepsilon
s_n h^{1/2}/c' \bigr)^{1/2}
\\
&&\qquad \le \bigl(c'h^{-1/2} \bigr)E \bigl\{
\epsilon^4 \bigr\}^{1/2} \bigl(\varepsilon^4
s_n^4 h^2 \bigr)^{-1/2}E \bigl\{
\varepsilon^4 \bigr\}^{1/2}
\\
&&\qquad =\frac{(c')^2E\{\varepsilon^4\}}{\varepsilon^2s_n^2h^2}\rightarrow0,
\end{eqnarray*}
where the last limit follows by $s_n^2\asymp n$ and the assumption
$nh^2\rightarrow\infty$.
Then as $n$ approaches infinity,
\begin{eqnarray*}
&&\frac{1}{s_n^2}\sum_{i=1}^nE \bigl
\{\bigl| \epsilon_i \bigl(x^T H_{U_i}+h^{1/2}
T_{U_i}(z_0) \bigr)\bigr|^2 I{\bigl|\epsilon_i
\bigl(x^T H_{U_i}+h^{1/2} T_{U_i}(z_0)
\bigr)\bigr|\ge\varepsilon s_n} \bigr\}
\\
&&\qquad =\frac{n}{s_n^2}E \bigl\{\bigl|\epsilon \bigl(x^T
H_U+h^{1/2} T_U(z_0)
\bigr)\bigr|^2 I{\bigl|\epsilon \bigl(x^T H_U+h^{1/2}
T_U(z_0) \bigr)\bigr|\ge\varepsilon s_n} \bigr\}
\rightarrow0.
\end{eqnarray*}
So Lindeberg's condition holds. The desired result follows immediately
by central limit theorem.
This completes the proof.
\end{pf*}

\subsection{Proof of Theorem \texorpdfstring{\protect\ref{mainthm3}}{3.1}}\label{prooflemma4}

The proof of Theorem~\ref{mainthm3} directly follows Theorem~\ref{mainthm2}
and the following lemma.

\begin{lemma}\label{lemma4} Suppose that there exists $b\in(1/(2m),1]$
such that $G_k$ satisfies~(\ref{smoothG}).
Then we have, for any $z_0\in\mathbb I$, $h^{1/2}A(z_0)=o(1)$,
$h^{1/2}(W_\lambda A)(z_0)=o(1)$. Furthermore, if $n^{1/2}h^{m(1+b)}=o(1)$,
then as $n\rightarrow\infty$,
\begin{equation}
\label{alimit}
\pmatrix{
\sqrt{n} \bigl(
\theta_0^*-\theta_0 \bigr)
\vspace*{3pt}\cr
\sqrt{nh} \bigl(g_0^*(z_0)-g_0(z_0)+(W_\lambda
g_0) (z_0) \bigr)}
\longrightarrow0.
\end{equation}
\end{lemma}

\begin{pf}
We will show (\ref{alimit}) in three steps:
\begin{longlist}[(iii)]
\item[(i)]  Show $\|V(G,W_\lambda g_0)\|_{l_2}=o(n^{-1/2})$. By (\ref{PropRPK}),
\[
V(G_k,W_\lambda g_0)=\sum
_{\mu\in\mathbb{Z}}V(G_k,h_\mu )V(g_0,h_\mu)
\frac
{\lambda\gamma_\mu}{1+\lambda\gamma_\mu},
\]
for any $k=1,\ldots,p$. Then by Cauchy's inequality, we have
\begin{eqnarray*}
&&\bigl|V(G_k,W_\lambda g_0)\bigr|^2
\\
&&\qquad \le\sum_{\mu}\bigl|V(G_k,h_\mu)\bigr|^2
\frac{\lambda\gamma_\mu
}{1+\lambda\gamma
_\mu} \sum_{\mu}\bigl|V(g_0,h_\mu)\bigr|^2
\frac{\lambda\gamma_\mu}{1+\lambda
\gamma_\mu
}
\\
&&\qquad\le \mathrm{const}\cdot\lambda\sum_{\mu}\bigl|V(G_k,h_\mu)\bigr|^2
\frac
{\lambda
\gamma_\mu}{1+\lambda\gamma_\mu}
\\
&&\qquad=\mathrm{const}\cdot\lambda\sum_{\mu}\bigl|V(G_k,h_\mu)\bigr|^2
\gamma _\mu ^b \biggl(\frac{\lambda\gamma_\mu^{1-b}}{1+\lambda\gamma_\mu
} \biggr)
\\
&&\qquad\le \mathrm{const}\cdot\lambda^{1+b}.
\end{eqnarray*}
Thus, when $n^{1/2}\lambda^{(1+b)/2}=n^{1/2}h^{m(1+b)}=o(1)$,
$\|V(G,W_\lambda g_0)\|_{l_2}=o(n^{-1/2})$.

\item[(ii)]  Show $\|A_k\|_{\sup}=O(1)$, for any $k=1,\ldots,p$. Note for any
$z\in\mathbb{I}$, by (\ref{inter1}),
\begin{eqnarray*}
A_k(z) &=& \langle A_k, K_z
\rangle_1 =V(G_k,K_z) \\
&=& \sum
_{\mu\in\mathbb{N}}\frac{V(G_k,h_\mu)}{1+\lambda\gamma_\mu} h_\mu(z).
\end{eqnarray*}
By boundedness of $h_\nu$s (Assumption~\ref{A.2}) and by
Cauchy's inequality, uniformly for $z\in\mathbb{I}$,
\begin{eqnarray*}
\bigl|A_k(z)\bigr|^2&\le&\sum_{\mu}\bigl|V(G_k,h_\mu)\bigr|^2(1+
\gamma_\mu)^b \bigl|h_\mu(z)\bigr|^2 \cdot
\sum_{\mu} \frac{1}{(1+\gamma_\mu)^b(1+\lambda\gamma_\mu
)^2}
\\
&=&O \biggl(\sum_{\mu} \frac{1}{(1+\gamma_\mu)^b}
\biggr)=O(1),
\end{eqnarray*}
where the last equality follows by $\gamma_\mu\asymp\mu^{2m}$ and $2mb>1$.
This shows\break $\|A_k\|_{\sup}=O(1)$, implying $h^{1/2}A(z_0)=o(1)$.
By (\ref{wlamA}), $(W_\lambda A)(z)=A(z)-\sum_\mu\frac{V(G,h_\mu
)}{(1+\lambda\gamma_\mu)^2}h_\mu(z)$.
Using the above derivations we can show that uniformly for $z\in
\mathbb{I}$,
$|\sum_\mu\frac{V(G,h_\mu)}{(1+\lambda\gamma_\mu)^2}h_\mu
(z)|^2=O
(\sum_{\mu} \frac{1}{(1+\gamma_\mu)^b} )=O(1)$,\vspace*{1pt}
implying\break $h^{1/2}(W_\lambda A)(z_0)=o(1)$.

\item[(iii)]  By (i) and (ii), (\ref{alimit}) follows by, as
$n\rightarrow\infty$,
\begin{eqnarray*}
&&\pmatrix{
n^{1/2} \bigl(
\theta_0^*-\theta_0 \bigr)
\vspace*{3pt}\cr
(nh)^{1/2} \bigl(g_0^*(z)-g_0(z)+(W_\lambda
g_0) (z) \bigr) }
\\
&&\qquad =\pmatrix{
n^{1/2}(\Omega+
\Sigma_\lambda)^{-1}V(G,W_\lambda g_0)
\vspace*{3pt}\cr
-(nh)^{1/2}V \bigl(G^T,W_\lambda g_0
\bigr) (\Omega+\Sigma_\lambda)^{-1}A(z)}
  \rightarrow0.
\end{eqnarray*}
\end{longlist}
\upqed\end{pf}
%

\subsection{Proof of Theorem \texorpdfstring{\protect\ref{LRTlimdist}}{4.4}}\label{secproofjointlocaltest}
For notational convenience, denote
$\widehat{f}=\widehat{f}_{n,\lambda}$,
$\widehat{f}^0=\widehat{f}_{n,\lambda}^{H_0}$, the constrained estimate
of $f$ under $H_0$, and
$f=\widehat{f}^0-\widehat{f}=(\theta,g)$. By Assumptions\vspace*{1pt}
\ref{A.5} and \ref{A.6}, with large probability,
$\|f\|\le r_n$, where $r_n=M ((nh)^{-1/2}+h^m)$ for some large
$M$. By Assumption~\ref{A.1}(a), for some large
constant $C>0$, the event $B_n\equiv B_{n1}\cap B_{n2}$ has large
probability, where $B_{n1}=\{\max_{1\le i\le
n}\sup_{a\in\mathcal{I}}|\ddot{\ell}_a(Y_i;a)|\le
C\log{n}\}$ and $B_{n2}=\{\max_{1\le i\le
n}\sup_{a\in\mathcal{I}}|\ell'''_a(Y_i;\break  a)|\le C\log{n}\}$.
Let $a_n$ be defined as in (\ref{bahareps}).

By Taylor's expansion,
\begin{eqnarray}
 \quad\mathrm{LRT}_{n,\lambda} &=&\ell_{n,\lambda} \bigl(
\widehat{f}^0 \bigr)-\ell_{n,\lambda}(\widehat {f})
\nonumber
\\
&=&S_{n,\lambda}(\widehat{f}) f+\int_0^1
\!\!\int_0^1 s\mathit{DS}_{n,\lambda
} \bigl(
\widehat{f}+ss'f \bigr)ff\,ds\,ds'
\nonumber
\\
\label{LocalLRTeq1}
&=&\int_0^1\!\!\int_0^1
s\mathit{DS}_{n,\lambda} \bigl(\widehat {f}+ss'f \bigr)ff\,da\,ds'
\\
\nonumber
&=&\int_0^1\!\!\int_0^1
s \bigl\{\mathit{DS}_{n,\lambda} \bigl(\widehat {f}+ss'f
\bigr)ff-\mathit{DS}_{n,\lambda
}(f_0)ff \bigr\}\,ds\,ds'
\nonumber
\\
&&{}+\frac{1}{2} \bigl(\mathit{DS}_{n,\lambda}(f_0)ff-E \bigl
\{\mathit{DS}_{n,\lambda}(f_0)ff \bigr\} \bigr)
+\frac{1}{2}E \bigl\{\mathit{DS}_{n,\lambda}(f_0)ff \bigr\}.\hspace*{-12pt}\nonumber
\end{eqnarray}
Denote\vspace*{1pt} the above three sums by $I_1$, $I_2$ and $I_3$. Next we
will study the asymptotic behavior of these sums. Denote
$\tilde{f}=\widehat{f}+ss'f-f_0=(\tilde{\theta},\tilde{g})$,
for any $0\le s,s'\le1$. So $\|\tilde{f}\|=O_{P}(r_n)$.

By calculations of the Frech\'{e}t derivatives, we have
\begin{eqnarray*}
&&\!\!\mathit{DS}_{n,\lambda} \bigl(\widehat{f}+ss'f \bigr)ff
\\
&& \!\!\qquad =\mathit{DS}_{n,\lambda}(\tilde{f}+f_0)ff
\\
&&\!\!\qquad =\frac{1}{n}\sum_{i=1}^n \ddot{
\ell}_a \bigl(Y_i; X_i^T
\theta_0+g_0(Z_i)+X_i^T
\tilde{\theta}+\tilde{g}(Z_i) \bigr) \bigl(X_i^T
\theta+g(Z_i) \bigr)^2-\langle P_\lambda f,f
\rangle,
\end{eqnarray*}
and
\[
\mathit{DS}_{n,\lambda}(f_0)ff=\frac{1}{n}\sum
_{i=1}^n \ddot{\ell}_a
\bigl(Y_i; X_i^T\theta_0+g_0(Z_i)
\bigr) \bigl(X_i^T\theta+g(Z_i)
\bigr)^2-\langle P_\lambda f,f\rangle.
\]
On $B_n$,
\begin{eqnarray}
&&\bigl|\mathit{DS}_{n,\lambda} \bigl(\widehat{f}+ss'f
\bigr)ff-\mathit{DS}_{n,\lambda
}(f_0)ff\bigr|
\nonumber
\\
&&\qquad\le\frac{1}{n}C(\log{n}) \|\tilde{f}\|_{\sup}\sum
_{i=1}^n \bigl(X_i^T
\theta+g(Z_i) \bigr)^2
\nonumber
\\
&&\qquad =C(\log{n})\|\tilde{f}\|_{\sup} \Biggl\langle\frac{1}{n}\sum
_{i=1}^n \bigl(X_i^T
\theta+g(Z_i) \bigr) R_{U_i},f \Biggr\rangle
\nonumber
\\[-8pt]
\label{LocalLRTeq2}
\\[-8pt]
\nonumber
&&\qquad =C(\log{n})\|\tilde{f}\|_{\sup} \Biggl\langle\frac{1}{n}\sum
_{i=1}^n \bigl(X_i^T
\theta+g(Z_i) \bigr) R_{U_i}\nonumber\\
&&\hspace*{40pt}\qquad\qquad\qquad{}-E_T \bigl\{
\bigl(X^T\theta+g(Z) \bigr)R_U \bigr\},f \Biggr\rangle
\nonumber
\\
&&\qquad\quad{}+C(\log{n})\|\tilde{f}\|_{\sup}E_T \bigl\{
\bigl(X^T\theta+g(Z) \bigr)^2 \bigr\}.\nonumber
\end{eqnarray}
Now we study
$\frac{1}{n}\|\sum_{i=1}^n(X_i^T\theta+g(Z_i))
R_{U_i}-E_T\{(X^T\theta+g(Z))R_U\}\|$. Let $d_n=c_m h^{-1/2}
r_n$ and $\bar{f}= d_n^{-1} f/2=(d_n^{-2}\theta/2,d_n^{-1}g/2)\equiv
(\bar{\theta},\bar{g})$. Consider\vspace*{1pt}
$\psi(T;f)=X^T\theta+g(Z)$ and $\psi_n(T;\bar{f})=(1/2)c_m^{-1}
h^{1/2}d_n^{-1}\psi(T;2d_n \bar{f})$. It is easy to see that $\psi
_n(T;\bar{f})$, as a function of $\bar{f}$, satisfies the
Lipschitz continuity condition (S.6) in the online supplementary.

Since $h=o(1)$ and $nh^2\rightarrow\infty$, $d_n=o(1)$. Then by
Lemma~\ref{lemma0}, on $B_n$, $\|\bar{f}\|_{\sup}\le1/2$,
which implies that for any $(x,z)\in\mathcal{U}$,
$|x^T\bar{\theta}+\bar{g}(z)|\le1/2$. Letting $x$
approach zero, one gets that $|\bar{g}(z)|\le1/2$, and
thus, $\|\bar{g}\|_{\sup}\le1/2$, which further implies
that $|x^T\bar{\theta}|\le
\|\bar{g}\|_{\sup}+\|\bar{f}\|_{\sup}\le1$ for any
$x\in\mathbb{I}^p$. Also note that
\begin{eqnarray*}
J(\bar{g},\bar{g})&=& d_n^{-2}\lambda^{-1}
\bigl( \lambda J(g,g) \bigr)/4
\\
&\le& d_n^{-2}\lambda^{-1} \|f\|^2/4
\\
&\le& d_n^{-2}\lambda^{-1} r_n^2/4
\\
&<&c_m^{-2} h\lambda^{-1}.
\end{eqnarray*}
Thus, when event $B_n$ holds, $\bar{f}$ is an element in $\mathcal{G}$.
Then by Lemma S.3 (in the supplementary material \cite{CS15}), with
large probability
\begin{eqnarray}
&&\Biggl\|\frac{1}{n}\sum_{i=1}^n
\bigl[ \bigl(X_i^T\theta+g(Z_i)
\bigr)R_{U_i}-E_T \bigl\{ \bigl(X^T\theta +g(Z)
\bigr)R_U \bigr\} \bigr]\Biggr\|
\nonumber
\\
\label{eqastep}
&&\qquad =\frac{c_mh^{-1/2}d_n}{n}\Biggl\|\sum_{i=1}^n
\bigl[\psi_n(T_i;\bar {f})R_{U_i}-E_T
\bigl\{\psi_n(T; \bar{f})R_U \bigr\} \bigr]\Biggr\|
\\
&&\qquad =O_{P} \bigl(a_n' \bigr),\nonumber
\end{eqnarray}
where $a_n'=n^{-1/2}((nh)^{-1/2}+h^m)
h^{-(6m-1)/(4m)}(\log\log{n})^{1/2}$.
So by $a_n'=o(r_n)$,
\begin{eqnarray}
&&\bigl|\mathit{DS}_{n,\lambda} \bigl(\widehat{f}+ss'f
\bigr)ff-\mathit{DS}_{n,\lambda
}(f_0)ff\bigr|
\nonumber
\\
&&\qquad = \|\tilde{f}\|_{\sup} \bigl(O_{P} \bigl(a_n'
r_n\log{n} \bigr)+O_{P} \bigl(r_n^2
\log {n} \bigr) \bigr)
\nonumber
\\[-8pt]
\label{eq1astep}
\\[-8pt]
\nonumber
&&\qquad =h^{-1/2} r_n O_{P} \bigl(r_n^2
\log{n} \bigr)
\nonumber
\\
&&\qquad =O_{P} \bigl(r_n^3 h^{-1/2}\log{n}
\bigr).\nonumber
\end{eqnarray}
Thus $|I_1|=O_{P}(r_n^3 h^{-1/2}\log{n})$.

Next we approximate $I_2$.
Define $\psi(T;f)=\ddot{\ell}_a(Y; X^T\theta_0+g_0(Z))(X^T\theta+g(Z))$.
Then by calculation of the Fr\'{e}chet derivative (Section~\ref
{subsection2.3}),
\begin{eqnarray*}
&& \mathit{DS}_{n,\lambda}(f_0)ff-E \bigl\{\mathit{DS}_{n,\lambda}(f_0)ff
\bigr\}\\
&& \qquad = \Biggl\langle\frac{1}{n}\sum_{i=1}^n
\bigl[\psi(T_i;f)R_{U_i}-E_T \bigl\{\psi
(T;f)R_U \bigr\} \bigr],f \Biggr\rangle.
\end{eqnarray*}
Thus $2|I_2|\le\frac{1}{n}\|\sum_{i=1}^n[\psi(T_i;f)R_{U_i}-E_T\{
\psi
(T;f)R_U\}]\|\cdot\|f\|$.
So it is sufficient to approximate\vspace*{1pt} $\|\sum_{i=1}^n[\psi
(T_i;f)R_{U_i}-E_T\{\psi(T;f)R_U\}]\|$.
Let $\tilde{\psi}_n(T;\bar{f})=\break (1/2)C^{-1} c_m^{-1} (\log{n})^{-1}h^{1/2}
d_n^{-1}\psi(T;2d_n \bar{f})$ and $\psi_n(T_i;\bar{f})=\tilde{\psi
}_n(T_i;\bar{f})
I_{A_i}$, where $\bar{f}=d_n^{-1}f/2$ and $A_i=\{\sup_{a\in\mathcal{I}}
|\ddot{\ell}_a(Y_i;a)|\le C\log{n}\}$ for $i=1,\ldots,n$. By similar
derivations as the ones below (\ref{LocalLRTeq2}),
it can be shown that
on $B_n$, $\bar{f}\in\mathcal{G}$. Observe
that $B_n$ implies $\bigcap_i A_i$.
A direct examination shows that $\psi_n$ satisfies (S.6).
By Lemma S.3, with large probability,
\begin{eqnarray}
&& \Biggl\|\sum_{i=1}^n \bigl[
\psi_n(T_i;\bar{f}) R_{U_i}-E_T
\bigl\{\psi_n(T;\bar{f})R_U \bigr\} \bigr]\Biggr\|
\nonumber
\\[-7pt]
\label{eq1LRTdistappendix}
\\[-7pt]
\nonumber
&& \qquad \le
\bigl(n^{1/2}h^{-(2m-1)/(4m)}+1 \bigr) (5\log\log{n})^{1/2}.
\end{eqnarray}

On the other hand, by Chebyshev's inequality
\begin{eqnarray*}
P \bigl(A_i^c \bigr)&=& \exp \bigl(-(C/C_0)
\log{n} \bigr) E \Bigl\{\exp \Bigl(\sup_{a\in\mathcal{I}}\bigl|\ddot{
\ell}_a(Y_i;a)\bigr|/C_0 \Bigr) \Bigr\}\\
&\le&
C_1n^{-C/C_0}.
\end{eqnarray*}
Since $h=o(1)$ and $nh^2\rightarrow\infty$, we may choose $C$ to
be large so that $(\log{n})^{-1} \times\break n^{-C/(2C_0)}=o( a_n' h^{1/2}
d_n^{-1})$,\vspace*{1pt} where
\[
a_n'=n^{-1/2} \bigl((nh)^{-1/2}+h^m
\bigr) h^{-(6m-1)/(4m)}(\log\log{n})^{1/2}.
\]
By (\ref{A1aeq1}), which implies
$E\{\sup_{a\in\mathcal{I}}|\ddot{\ell}_a(Y_i;a)| |
U_i\}\le2C_1C_0^2$, we have, on $B_n$,
$E_T\{|\psi(T;2d_n\bar{f})|^2\}\le2C_1C_0^2 d_n^2$.
So when $n$ is large, on $B_n$, by Chebyshev's inequality,\vspace*{1pt}
\begin{eqnarray}
&&\bigl\|E_T \bigl\{\psi_n(T_i;
\bar{f})R_{U_i} \bigr\}-E_T \bigl\{\tilde{
\psi}_n(T_i; \bar{f}) R_{U_i} \bigr\}\bigr\|
\nonumber
\\
&&\qquad=\bigl\|E_T \bigl\{\tilde{\psi}_n(T_i;\bar{f})
R_{U_i} \cdot I_{A_i^c} \bigr\}\bigr\|
\nonumber
\\
\label{eq2LRTdistappendix}
&&\qquad \le (1/2)C^{-1} (\log{n})^{-1} d_n^{-1}
\bigl(E_T \bigl\{\bigl|\psi(T;2d_n \bar {f})\bigr|^2
\bigr\} \bigr)^{1/2} P \bigl(A_i^c
\bigr)^{1/2}
\\
&&\qquad \le (1/2)2^{1/2}C^{-1}C_0C_1 (
\log{n})^{-1}n^{-C/(2C_0)}
\nonumber
\\
&&\qquad = o \bigl(a_n' h^{1/2} d_n^{-1}
\bigr).\nonumber
\end{eqnarray}
Therefore, by (\ref{eq1LRTdistappendix}) and (\ref
{eq2LRTdistappendix}), on $B_n$ with large\vspace*{0pt} probability,
\begin{eqnarray}
&&\frac{1}{n}\Biggl\|\sum_{i=1}^n
\bigl[\psi(T_i;f)R_{U_i}-E_T \bigl\{
\psi(T;f)R_U \bigr\} \bigr]\Biggr\|
\nonumber
\\
\nonumber
&&\qquad =\frac{2C c_m (\log{n}) h^{-1/2}d_n}{n}\Biggl\|\sum_{i=1}^n
\bigl[ \tilde{\psi }_n(T_i; \bar{f})R_{U_i}-E_T
\bigl\{\tilde{\psi}_n(T;\bar{f})R_U \bigr\} \bigr]\Biggr\|
\\
\label{LocalLRTeq3}
&&\qquad\le\frac{2C c_m (\log{n}) h^{-1/2}d_n}{n} \\
&&\qquad\quad{}\times\Biggl(\Biggl\|\sum_{i=1}^n
\bigl[\psi _n(T_i;\bar{f})R_{U_i}-E_T
\bigl\{\psi_n(T;\bar{f})R_U \bigr\} \bigr]\Biggr\|\nonumber
\\
&&\hspace*{18pt}\qquad\quad{}+n\bigl\|E_T \bigl\{\psi_n(T_i;
\bar{f})R_{U_i} \bigr\}-E_T \bigl\{\tilde{\psi
}_n(T_i;\bar {f}) R_{U_i} \bigr\}\bigr\| \Biggr)
\nonumber
\\
&&\qquad \le\frac{2C c_m (\log{n}) h^{-1/2}d_n}{n}
\nonumber
\\
&&\quad\qquad{}\times \bigl[ \bigl(n^{1/2}h^{-(2m-1)/(4m)}+1 \bigr) (5\log
\log{n})^{1/2}+o \bigl(n a_n' h^{1/2}
d_n^{-1} \bigr) \bigr]
\nonumber
\\
&&\qquad \le C'a_n' \log{n},
\nonumber
\end{eqnarray}
for some large constant $C'>0$.
Thus $|I_2|=O_{P}(a_n'r_n\log{n})$.

Note that $I_3=-\|f\|^2/2$.
Therefore,
\begin{eqnarray*}
-2n\cdot \mathrm{LRT}_{n,\lambda}&=&n \bigl\|\widehat{f}^0-\widehat{f}
\bigr\|^2+O_{P} \bigl(n r_n a_n'
\log{n}+nr_n^3h^{-1/2}\log{n} \bigr)
\\
&=&n \bigl\|\widehat{f}^0-\widehat{f}\bigr\|^2+O_{P}
\bigl(n r_n a_n\log{n}+nr_n^3h^{-1/2}
\log{n} \bigr).
\end{eqnarray*}
By $r_n^2h^{-1/2}=o(a_n)$ and
$nr_na_n=o((\log{n})^{-1})$, we have that $O_{P}(n r_n
a_n\log{n}+nr_n^3h^{-1/2}\log{n})=o_{P}(1)$. This shows $-2n\cdot
\mathrm{LRT}_{n,\lambda}=n
\|\widehat{f}^0-\widehat{f}\|^2+o_P(1)$. So we only focus on
$n \|\widehat{f}^0-\widehat{f}\|^2$. By Theorems
\ref{mainthm1} and \ref{LRTmainthm1},
\begin{equation}
\label{approxeq0}
 \quad n^{1/2}\bigl\|\widehat{f}^0-
\widehat{f}-S_{n,\lambda
}^0 \bigl(f_0^0
\bigr)+S_{n,\lambda
}(f_0) \bigr\|=O_{P} \bigl(n^{1/2}a_n
\log{n} \bigr)=o_{P}(1),\hspace*{-6pt}
\end{equation}
so we just have to focus on
$n^{1/2}\{S_{n,\lambda}^0(f_0^0)-S_{n,\lambda}(f_0)\}$.
Recall that under $H_0$, $f_0^0=(\theta_0^0,g_0^0)\in\mathcal{H}_0$, so
\begin{eqnarray*}
S_{n,\lambda}^0 \bigl(f_0^0 \bigr) &=&
\frac{1}{n}\sum_{i=1}^n\dot{
\ell}_a \bigl(Y_i; X_i^T\theta
_0^0+g_0^0(Z_i)+X_i^T
\theta^\dag+w^\dag \bigr)R_{U_i}^0-P^0_\lambda
f_0^0
\\
&=&\frac{1}{n}\sum_{i=1}^n\dot{
\ell}_a \bigl(Y_i; X_i^T\theta
_0+g_0(Z_i) \bigr)R_{U_i}^0-P^0_\lambda
f_0^0
=\frac{1}{n}\sum_{i=1}^n
\epsilon_i R_{U_i}^0-P^0_\lambda
f_0^0,
\end{eqnarray*}
where $\epsilon_i=\dot{\ell}_a(Y_i;X_i^T\theta_0+g_0(Z_i))$, $R_U^0$
and $P^0_\lambda f_0^0$ are defined in Section~\ref{seclrt},
and
\[
S_{n,\lambda}(f_0)=\frac{1}{n}\sum
_{i=1}^n \epsilon_i
R_{U_i}-P_\lambda f_0.
\]
Consequently,
\begin{eqnarray*}
&&S_{n,\lambda}^0 \bigl(f_0^0
\bigr)-S_{n,\lambda}(f_0)
\\
&&\qquad =\frac{1}{n}\sum_{i=1}^n
\epsilon_i \bigl(R_{U_i}^0-R_{U_i}
\bigr)- \bigl(P^0_\lambda f_0^0-P_\lambda
f_0 \bigr)
\\
&&\qquad =-\frac{1}{n}\sum_{i=1}^n
\epsilon_i \Biggl(\sum_{j=1}^k
\rho _{U_i,j}R_{q_j,W_j} \Biggr)+ \Biggl(\sum
_{j=1}^k\zeta_jR_{q_j,W_j}
\Biggr)
\\
&&\qquad =-\frac{1}{n}\sum_{i=1}^n
\epsilon_i \bigl(H(Q,W)\rho_{U_i}, \bigl(Q^TK_{z_0}-A^TH(Q,W)
\bigr)\rho_{U_i} \bigr)
\\
&&\qquad\quad{}+ \bigl(H(Q,W)\zeta, \bigl(Q^TK_{z_0}-A^TH(Q,W)
\bigr)\zeta \bigr)
\\
&&\qquad=(\xi, \beta)+ \bigl(H(Q,W)\zeta, \bigl(Q^TK_{z_0}-A^TH(Q,W)
\bigr)\zeta \bigr),
\end{eqnarray*}
where
\[
\beta=-\delta K_{z_0}-A^T\xi, \xi=-(1/n)\sum
_{i=1}^n\epsilon_i H(Q,W)
\rho_{U_i}
\]
and
\[
\delta=(1/n)\sum_{i=1}^n
\epsilon_iQ^T\rho_{U_i}.
\]
Therefore,
\begin{eqnarray*}
&&\bigl\|S_{n,\lambda}^0 \bigl(f_0^0
\bigr)-S_{n,\lambda}(f_0)\bigr\|^2
\\
&&\qquad =\bigl\|(\xi,\beta)\bigr\|^2+2 \bigl\langle(\xi,\beta), \bigl(H(Q,W)\zeta,
\bigl(Q^TK_{z_0}-A^TH(Q,W) \bigr)\zeta \bigr)
\bigr\rangle
\\
&&\qquad\quad{}+\bigl\| \bigl(H(Q,W)\zeta, \bigl(Q^TK_{z_0}-A^TH(Q,W)
\bigr)\zeta \bigr)\bigr\|^2.
\end{eqnarray*}
We next evaluate the three items on the right-hand side of the above
equation. Denote
$\Sigma_\lambda=E_U\{I(U)(G(Z)-A(Z))(G(Z)-A(Z))^T\}$. Note\break 
$E_Z\{B(Z) (G(Z)-A(Z))K_{z_0}(Z)\}=V(G,K_{z_0})-V(A,K_{z_0})=\langle
A,K_{z_0}\rangle_1-\break V(A,K_{z_0})
=\langle W_\lambda A,  K_{z_0}\rangle_1=(W_\lambda A)(z_0)$.
First,
\begin{eqnarray}
&&\bigl\|(\xi,\beta)\bigr\|^2
\nonumber
\\
&&\qquad =E_U \bigl\{I(U) \bigl(X^T\xi+\beta(Z)
\bigr)^2 \bigr\}+\lambda J(\beta,\beta)
\nonumber
\\
&&\qquad =E_U \bigl\{I(U) \bigl[ \bigl(X-A(Z) \bigr)^T\xi-
\delta K_{z_0}(Z) \bigr]^2 \bigr\}+\lambda J(\beta ,\beta )
\nonumber
\\
&&\qquad =\xi^T E_U \bigl\{I(U) \bigl(X-A(Z) \bigr)
\bigl(X-A(Z) \bigr)^T \bigr\}\xi
\nonumber
\\
&&\qquad\quad{}-2\xi^T E_U \bigl\{I(U) \bigl(X-A(Z)
\bigr)K_{z_0}(Z) \bigr\}\delta
\nonumber
\\
&&\label{evalitem1}\qquad\quad{}+\delta^2 E_Z{B(Z)\bigl|K_{z_0}(Z)\bigr|^2}
+ \bigl\langle W_\lambda \bigl(\delta K_{z_0}+A^T\xi
\bigr), \delta K_{z_0}+A^T\xi \bigr\rangle _1
\\
&&\qquad=\xi^T(\Omega+\Sigma_\lambda)\xi-2\xi^T
E_Z \bigl\{ B(Z) \bigl(G(Z)-A(Z) \bigr)K_{z_0}(Z) \bigr\}
\delta
\nonumber
\\
&&\qquad\quad{}+\delta^2 V(K_{z_0},K_{z_0})+\delta^2\langle W_\lambda K_{z_0},
K_{z_0}\rangle_1\nonumber\\
&&\qquad\quad{}+2\delta\xi ^T\langle
W_\lambda A, K_{z_0}\rangle_1 +\xi^T
\bigl\langle W_\lambda A, A^T \bigr\rangle_1\xi
\nonumber
\\
&&\qquad=\xi^T\Gamma_\lambda\xi-2\xi^T
(W_\lambda A) (z_0)\delta+\delta^2
K(z_0,z_0)+2\delta\xi^T (W_\lambda A)
(z_0)
\nonumber
\\
&&\qquad=\xi^T\Gamma_\lambda\xi+\delta^2
K(z_0,z_0),
\nonumber
\end{eqnarray}
where $\Gamma_\lambda=\Omega+\Sigma_\lambda+\langle W_\lambda A,
A^T\rangle_1$ and $\Sigma_\lambda=E_Z\{B(Z)(G(Z)-A(Z))(G(Z)-A(Z))^T\}$.
Second,
\begin{eqnarray}
&& \bigl\langle(\xi,\beta), \bigl(H(Q,W)\zeta,
\bigl(Q^TK_{z_0}-A^TH(Q,W) \bigr)\zeta \bigr)
\bigr\rangle
\nonumber
\\
&& \qquad =E_U \bigl\{I(U) \bigl[ \bigl(X-A(Z) \bigr)^T\xi-
\delta K_{z_0}(Z) \bigr]
\nonumber
\\
&&\qquad\hspace*{18pt}\quad{}\times \bigl[ \bigl(X-A(Z) \bigr)^T H(Q,W)\zeta+Q^T
\zeta K_{z_0}(Z) \bigr] \bigr\}
\nonumber
\\
&&\qquad\quad{}+ \bigl\langle W_\lambda\beta,Q^T\zeta
K_{z_0}-A^T H(Q,W)\zeta \bigr\rangle _1
\nonumber
\\
&&\qquad=\xi^T E_U \bigl\{I(U) \bigl(X-A(Z) \bigr)
\bigl(X-A(Z) \bigr)^T \bigr\}H(Q,W)\zeta
\nonumber
\\
\label{evalitem2}
&&\qquad\quad{}+\xi^T E_U \bigl\{I(U) \bigl(X-A(Z)
\bigr)K_{z_0}(Z) \bigr\}Q^T\zeta
\\
&&\qquad\quad{}-\delta E_U \bigl\{I(U)K_{z_0}(Z) \bigl(X-A(Z)
\bigr)^T \bigr\} H(Q,W)\zeta
\nonumber
\\
&&\qquad\quad{}-\delta Q^T\zeta V(K_{z_0},K_{z_0}) -\delta
Q^T\zeta\langle W_\lambda K_{z_0}, K_{z_0}
\rangle_1
\nonumber
\\
&&\qquad\quad{}+\delta \bigl(H(Q,W)\zeta \bigr)^T (W_\lambda A)
(z_0)-Q^T\zeta\xi^T(W_\lambda A)
(z_0)
\nonumber
\\
&&\qquad\quad{}+\xi^T \bigl\langle W_\lambda A, A^T \bigr
\rangle_1 H(Q,W)\zeta
\nonumber
\\
&&\qquad =\xi^T\Gamma_\lambda H(Q,W)\zeta-\delta Q^T
\zeta K(z_0,z_0).
\nonumber
\end{eqnarray}
Third, similar to the calculations in (\ref{evalitem1}) and (\ref
{evalitem2}), we have
\begin{eqnarray}
&& \bigl\langle \bigl(H(Q,W)\zeta, \bigl(Q^TK_{z_0}-A^TH(Q,W)
\bigr)\zeta \bigr),
 \nonumber\\
&&\hspace*{5pt}  \bigl(H(Q,W)\zeta,\bigl(Q^TK_{z_0}-A^TH(Q,W)
\bigr)\zeta \bigr) \bigr\rangle
\nonumber
\\
&&\label{evalitem3}\qquad =E_U \bigl\{I(U) \bigl[ \bigl(X-A(Z) \bigr)^TH(Q,W)
\zeta+Q^T\zeta K_{z_0}(Z) \bigr]^2 \bigr\}
\\
\nonumber
&&\qquad\quad{}+ \bigl\langle W_\lambda \bigl(Q^T\zeta
K_{z_0}-A^T H(Q,W)\zeta \bigr), Q^T\zeta
K_{z_0}-A^T H(Q,W) \zeta \bigr\rangle_1\hspace*{-12pt}
\\
&&\qquad =\zeta^T H(Q,W)^T\Gamma_\lambda H(Q,W)\zeta+
\bigl(Q^T\zeta \bigr)^2 K(z_0,z_0).\nonumber
\end{eqnarray}
It follows from (\ref{evalitem1}) to (\ref{evalitem3}) that
\begin{eqnarray}
&&\bigl\|S_{n,\lambda}^0 \bigl(f_0^0
\bigr)-S_{n,\lambda}(f_0)\bigr\|^2
\nonumber
\\
\label{importanteq}
&&\qquad = \bigl(\xi+H(Q,W)\zeta \bigr)^T\Gamma_\lambda \bigl(
\xi+H(Q,W)\zeta \bigr)+ \bigl(\delta -Q^T\zeta \bigr)^2
K(z_0,z_0)\hspace*{-18pt}
\\
&&\qquad =\pmatrix{\xi+H(Q,W)\zeta\vspace*{3pt}\cr
\delta-Q^T\zeta}^T \pmatrix{
\Gamma_\lambda& 0
\vspace*{3pt}\cr
0 & K(z_0,z_0) }
\pmatrix{
\xi+H(Q,W)\zeta \vspace*{3pt}\cr \delta-Q^T\zeta}.\nonumber
\end{eqnarray}
Next we find
the limiting distribution of $n\|S_{n,\lambda}^0(f_0^0)-S_{n,\lambda
}(f_0)\|^2$,
which leads to the limiting distribution of $-2n\cdot \mathrm{LRT}_{n,\lambda}$
in view of (\ref{approxeq0}).
By definition of $\xi$ and the expressions of $H(Q,W)$, $T(Q,W)$,
$\rho
_{U_i}$ and $\zeta$ in Section~\ref{seclrt}, we have
\begin{eqnarray*}
&&\xi+H(Q,W)\zeta
\\
&&\qquad =-\frac{1}{n}\sum_{i=1}^n
\epsilon_i H(Q,W) M_K^{-1}
\bigl(MH_{U_i}+QT_{U_i}(z_0) \bigr)
\\
&&\qquad\quad{}+ H(Q,W)M_K^{-1} \bigl(MH_{g_0}^*+QT_{g_0}^*(z_0)
\bigr)
\\
&&\qquad =H(Q,W)M_K^{-1}N \Biggl(-\frac{1}{n}\sum
_{i=1}^n\epsilon _i
\pmatrix{H_{U_i}  \vspace*{3pt}\cr T_{U_i}(z_0)}+ \pmatrix{H_{g_0}^*\vspace*{3pt}\cr
T_{g_0}^*(z_0)} \Biggr)
\\
&&\qquad=H(Q,W)M_K^{-1}N \pmatrix{
I_p&0
\vspace*{3pt}\cr
-A(z_0)^T&1 }
\\
&&\qquad\quad{}\times \Biggl(-\frac{1}{n}\sum_{i=1}^n
\epsilon_i \pmatrix{H_{U_i}\vspace*{3pt}\cr K_{z_0}(Z_i)}+
\pmatrix{H_{g_0}^*\vspace*{3pt}\cr (W_\lambda g_0) (z_0)}
\Biggr).
\end{eqnarray*}
On the other hand,
\begin{eqnarray*}
&&\delta-Q^T\zeta
\\
&& \qquad =Q^T M_K^{-1}N \Biggl(\frac{1}{n}\sum
_{i=1}^n\epsilon _i
\pmatrix{H_{U_i} \vspace*{3pt}\cr T_{U_i}(z_0)}-\pmatrix{H_{g_0}^*\vspace*{3pt}\cr
T_{g_0}^*(z_0)} \Biggr)
\\
&&\qquad =Q^T M_K^{-1}N
\pmatrix{
I_p&0
\vspace*{3pt}\cr
-A(z_0)^T&1}
\Biggl(\frac{1}{n}\sum_{i=1}^n
\epsilon_i \pmatrix{H_{U_i}\vspace*{3pt}\cr K_{z_0}(Z_i)}-
\pmatrix{H_{g_0}^*\vspace*{3pt}\cr (W_\lambda g_0) (z_0)}
\Biggr).
\end{eqnarray*}
Therefore,
\begin{eqnarray}
&&\pmatrix{\xi+H(Q,W)\zeta\vspace*{3pt}\cr\delta-Q^T\zeta}
\nonumber
\\
\label{acooleq}
&&\qquad =\pmatrix{H(Q,W)\vspace*{3pt}\cr-Q^T}M_K^{-1}N
\pmatrix{
I_p&0
\vspace*{3pt}\cr
-A(z_0)^T&1 }
\\
&&\qquad\quad{}\times \Biggl(-\frac{1}{n}\sum_{i=1}^n
\epsilon_i\pmatrix{H_{U_i}\vspace*{3pt}\cr K_{z_0}(Z_i)}+
\pmatrix{H_{g_0}^*\vspace*{2pt}\cr (W_\lambda g_0) (z_0)}
\Biggr).\nonumber
\end{eqnarray}
Define $\tilde{M}_K=
{H(Q,W)\choose  -Q^T}^T
{\hspace*{-9pt}\!\!\!\Gamma_\lambda \hspace*{9pt}\ \  \ 0\choose
0 \ \   \  K(z_0,z_0)}
{ H(Q,W)\choose  -Q^T
}$,
where\vspace*{1pt} we recall that
$\Gamma_\lambda=\Omega+\Sigma_\lambda+\langle
W_\lambda A,A^T\rangle_1$. Since\vspace*{1pt} for any $1\le j,k\le p$, $\langle
W_\lambda A_k,A_j\rangle_1=\lambda\sum_\nu
V(A_j, h_\nu)\times\break V(A_k,h_nu)\gamma_\nu=O(\lambda)=o(1)$, we have as
$\lambda\rightarrow0$, $\langle W_\lambda A,A^T\rangle_1\rightarrow
0$, a
$p\times p$ zero matrix. Define $\lambda_1$ as the maximum eigenvalue
of $\langle W_\lambda A,A^T\rangle_1$,
and $\lambda_2$ as the minimum eigenvalue of $\Omega+\Sigma_\lambda$.
Thus $\lambda_1=o(1)$. By equation (\ref{lemma3i}) in Lemma~\ref{lemma3},
$\lambda_2$ is asymptotically finitely upper bounded, and is lower
bounded from zero.
Note\vspace*{-2pt} that
\begin{eqnarray}
\tilde{M}_K-M_K
&=& \pmatrix{H(Q,W)
\vspace*{2pt}\cr
-Q^T }^T
\pmatrix{
\bigl\langle W_\lambda A,A^T
\bigr\rangle_1&0
\vspace*{2pt}\cr
0& 0 }
  \pmatrix{
H(Q,W)
\vspace*{2pt}\cr
-Q^T }
\nonumber
\\[-1pt]
&\le&\frac{\lambda_1}{\lambda_2}
\pmatrix{
H(Q,W)
\vspace*{2pt}\cr
-Q^T }^T\pmatrix{
\Omega+\Sigma_\lambda&0
\vspace*{2pt}\cr
0& 0 }
 \pmatrix{
H(Q,W)
\vspace*{2pt}\cr
-Q^T }
\nonumber
\\[-9pt]
\label{eqerror}
\\[-9pt]
\nonumber
&\le&\frac{\lambda_1}{\lambda_2}
\pmatrix{
H(Q,W)
\vspace*{2pt}\cr
-Q^T}^T
\pmatrix{
\Omega+\Sigma_\lambda&0
\vspace*{2pt}\cr
0& K(z_0,z_0)}
\pmatrix{
H(Q,W)
\vspace*{2pt}\cr
-Q^T }
\nonumber
\\[-1pt]
&=&\frac{\lambda_1}{\lambda_2}M_K.\nonumber
\end{eqnarray}
Define\vspace*{-2pt}
\begin{eqnarray*}
\Psi_\lambda
&=&\pmatrix{
(\Omega+\Sigma_\lambda)^{-1/2}&0
\vspace*{3pt}\cr
0&K(z_0,z_0)^{1/2} } \pmatrix{
I_p&-A(z_0)
\vspace*{3pt}\cr
0&1 }
  N^TM_K^{-1}
\tilde{M}_KM_K^{-1}N
\\
&&{}\times \pmatrix{
I_p&0
\vspace*{3pt}\cr
-A(z_0)^T&1 }
\pmatrix{
(\Omega+\Sigma_\lambda)^{-1/2}&0
\vspace*{3pt}\cr
0&K(z_0,z_0)^{1/2}}.
\end{eqnarray*}
Therefore, by (\ref{eqerror}),
\begin{eqnarray*}
0&\le&\Psi_\lambda-\Phi_\lambda
\\
&\le&\frac{\lambda_1}{\lambda_2} \pmatrix{
(\Omega+
\Sigma_\lambda)^{-1/2}&0
\vspace*{3pt}\cr
0&K(z_0,z_0)^{1/2} }
 \pmatrix{
I_p&-A(z_0)
\vspace*{3pt}\cr
0&1 }
  N^TM_K^{-1}N
\\
&&{}\times
\pmatrix{
I_p&0
\vspace*{3pt}\cr
-A(z_0)^T&1 }
  \pmatrix{
(\Omega+\Sigma_\lambda)^{-1/2}&0
\vspace*{3pt}\cr
0&K(z_0,z_0)^{1/2}}
  =\frac{\lambda_1}{\lambda_2}\Phi_\lambda.
\end{eqnarray*}
Since $\Phi_\lambda\le\textrm{trace}(\Phi_\lambda) I_{p+1}=k I_{p+1}$,
$\Psi_\lambda-\Phi_\lambda=o(1)I_{p+1}$. Thus, as $n\rightarrow
\infty$,
$\Psi_\lambda$ approaches $\Phi_0$.

Next\vspace*{1pt} we will complete the proof by demonstrating the asymptotic distribution.
It follows by Lemma~\ref{lemma4} that $n^{1/2}H_{g_0}^*=o(1)$.
Denote $N_U=(-H_U^T, K_Z(z_0)/K(z_0, z_0)^{1/2})^T$.
By Assumption~\ref{A.1}(c),
\begin{eqnarray*}
&&E \bigl\{\epsilon^2 N_U N_U^T
\bigr\}
\\
&&\qquad =E\lleft\{I(U) \pmatrix{
H_U H_U^T & -H_U
K_Z(z_0)/K(z_0,z_0)^{1/2}
\vspace*{3pt}\cr
-H_U^T K_Z(z_0)/K(z_0,z_0)^{1/2}
& |K_Z(z_0)|^2/K(z_0,z_0)
}
 \rright\}.
\end{eqnarray*}
To find the limit of this matrix, note that as
$\lambda\rightarrow0$, the following limits hold:
\begin{itemize}
\item by Lemma~\ref{lemma3},
\begin{eqnarray*}
&&E \bigl\{I(U)H_UH_U^T \bigr\}
\\[-2pt]
&&\qquad =(\Omega+\Sigma_\lambda)^{-1}E \bigl\{I(U) \bigl(X-A(Z)
\bigr) \bigl(X-A(Z) \bigr)^T \bigr\}(\Omega +\Sigma
_\lambda)^{-1}
\\[-2pt]
&&\qquad=(\Omega+\Sigma_\lambda)^{-1} \bigl(\Omega+E_Z
\{ B(Z) \bigl(G(Z)-A(Z) \bigr) \bigl(G(Z)-A(Z) \bigr)^T \bigr)
\\[-2pt]
&&\qquad\quad{}\times
(\Omega+\Sigma_\lambda)^{-1}
\\[-2pt]
&&\qquad \rightarrow \Omega^{-1};
\end{eqnarray*}
\item by $h^{1/2}(W_\lambda A)(z_0)\rightarrow0$ (see Lemma~\ref
{lemma4}) and
$hK(z_0,z_0)\rightarrow\sigma_{z_0}^2/c_0$ [by assumption (\ref
{LRTacondition})],
\begin{eqnarray*}
&&E \bigl\{I(U) H_U K_Z(z_0) \bigr
\}/K(z_0,z_0)^{1/2}
\\
&&\qquad =E \bigl\{I(U) (\Omega+\Sigma_\lambda)^{-1} \bigl(X-A(Z)
\bigr)K_Z(z_0) \bigr\} /K(z_0,z_0)^{1/2}
\\
&&\qquad =E \bigl\{B(Z) \bigl(G(Z)-A(Z) \bigr)K_Z(z_0) \bigr
\}/K(z_0,z_0)^{1/2}
\\
&&\qquad =(W_\lambda A) (z_0)/K(z_0,z_0)^{1/2}
\rightarrow0;
\end{eqnarray*}
\item by assumption, $E\{B(Z)|K_Z(z_0)|^2\}/K(z_0,z_0)\rightarrow
c_0$.
\end{itemize}
Thus, as $\lambda\rightarrow0$, $E\{\epsilon^2 N_U
N_U^T\}\rightarrow {
\Omega^{-1} \ \  0\choose
0 \ \ \ \ \ c_0}$. So as $n\rightarrow\infty$,
\begin{eqnarray}
&&\qquad n^{1/2}\pmatrix{
(\Omega+\Sigma_\lambda)^{1/2}&0
\vspace*{3pt}\cr
0&1 }
  \lleft(-\frac{1}{n}\sum
_{i=1}^n \epsilon_i \pmatrix{
H_{U_i}
\vspace*{3pt}\cr
\displaystyle\frac{K_{z_0}(Z_i)}{\sqrt{K(z_0,z_0)}}} + \pmatrix{
H_{g_0}^*
\vspace*{3pt}\cr
\displaystyle\frac{(W_\lambda g_0)(z_0)}{\sqrt{K(z_0,z_0)}}}
  \rright)\hspace*{-12pt}
\nonumber
\\[-8pt]
\label{acoollim}
\\[-8pt]
\nonumber
&&\qquad\qquad \stackrel{d} {\longrightarrow} \upsilon,
\end{eqnarray}
where $\upsilon\sim N\bigl(
{0\choose  c_{z_0}},
{I_p \ \  0\choose
0 \ \ c_0}
\bigr)$.
Therefore, it follows by (\ref{importanteq}), (\ref{acooleq}) and
(\ref{acoollim})
that, as $n\rightarrow\infty$, $n\|S_{n,\lambda
}^0(f_0^0)-S_{n,\lambda
}(f_0)\|^2\stackrel{d}{\longrightarrow} \upsilon^T \Phi_0 \upsilon$.
It immediately follows that $\|\widehat{f}^0-\widehat{f}\|=O_{P}(n^{-1/2})$.
Besides, when $n\rightarrow\infty$,
$-2n\cdot \mathrm{LRT}_{n,\lambda}
\stackrel{d}{\longrightarrow} \upsilon^T \Phi_0 \upsilon$.
\end{appendix}

\begin{supplement}[id=suppA]
\stitle{Supplement to
``Joint asymptotics for semi-nonparametric regression models with
partially linear structure''}
\slink[doi]{10.1214/15-AOS1313SUPP} 
\sdatatype{.pdf}
\sfilename{AOS1313\_supp.pdf}
\sdescription{Additional proofs are provided.}
\end{supplement}

\printaddresses
\end{document}